\documentstyle[11pt, amscd,amssymb,verbatim]{amsbook}
\newtheorem{thm}{Theorem}
\newtheorem{Cor}{Corollary}

\newtheorem{Prop}{Proposition}

\numberwithin{equation}{section}
\numberwithin{Cor}{Prop}
\numberwithin{Prop}{section}
\newtheorem{Rem}{Remark}
\numberwithin{Rem}{section}

\numberwithin{Def}{section}

\begin{document}

\title[Algebraicization]
{On the algebraicization of certain Stein manifolds. }
\author{R.  M. Aguilar}
\author{D. M. Burns, Jr.}
\address{Department of Mathematics \\ Tufts University \\ Medford, MA 02155 USA}
\email{aguilar@@umich.edu}
\address{Department of Mathematics \\ University of Michigan \\ Ann Arbor, MI 48109-1109, USA}
\email{dburns@@math.lsa.umich.edu}
\keywords{Riemannian manifold, CR manifold, Harmonic, 
Grauert tube} 

\begin{abstract}

To every real analytic Riemannian manifold $M$ there is associated a
complex structure on a neighborhood of the zero section in the real
tangent bundle of $M$. This structure can be uniquely specified in
several ways, and is referred to as a Grauert tube. We say that a
Grauert tube is entire if the complex structure can be extended to the
entire tangent bundle. We prove here that the complex manifold given
by an entire Grauert tube is, in a canonical way, an affine algebraic
variety. In the special case $M = S^2$, we show that any entire
Grauert tube associated to a metric (not necessarily round) on $M$
must be algebraically biholomorphic to the Grauert tube of the round
metric, that is, the non-singular quadric surface in ${\Bbb
C}^3$. (This second result has been discovered independently by Totaro.)

\end{abstract}
\maketitle

\section{{\bf Introduction}}
By a {\em Grauert tube} we mean 
a smooth complex manifold ${X}$ of complex dimension $n$
endowed with a  $C^5$ non-negative strictly plurisubharmonic
exhaustion 
function $\tau\colon {X}\to [0, R)$ whose 
square root is plurisubharmonic and 
satisfies the complex homogeneous Monge-Amp\`{e}re 
equation, 
\begin{equation}\label{MAESR}
\partial \Bar{\partial} (\sqrt{\tau})^n=0,
\end{equation}
on ${X}\setminus 
\{ \tau = 0\}$. 
That $\tau$ is an exhaustion means that 
for all $0\leq r <R$, $\{x\in {X}\mid 0\leq \tau(x) \leq r^2\}$ is 
compact.
If the maximal radius $R=\infty$ we call the tube 
{\em unbounded} or {\em 
entire}. The simplest example of such is the affine $2$-quadric 
$$
{\cal Q}=\{(z_1,z_2,z_3) 
\in {\mathbf C}^3 \mid  z_1^2 + z^2_2 + z^2_3=1\} \subset {\mathbf C}^3.
$$
with $$\sqrt{\tau} = \frac{1}{2} \cosh^{-1}(\|z\|^2).$$

By a classical theorem of Remmert, since $\tau$ is strictly
plurisubharmonic and ${X}$ is Stein, $X$ can be properly holomorphically
embedded in $\mathbf{C}^{N}$ for $N>>1$. The main result in this paper
is that an unbounded Monge-Amp\`{e}re exhaustion ensures the existence
of an {\em algebraic} embedding in ${\mathbf C}^N$, similar to the
simple example ${\cal Q}$ above. We show,
    
\begin{thm}
Any entire Grauert tube is an affine algebraic manifold.
\end{thm}

The functions 
giving the algebraic structure to ${X}$ can be described as the 
field of quotients of holomorphic functions $f$ on ${X}$
satisfying
$$\lim \sup_{r\to \infty}
\frac{1}{r}\int_{\{f(\tau) =r\}} \log_{+}\mid f \mid 
\tau^{-\frac{n}{2}} \ d\tau \wedge 
(\sqrt{-1} \partial\Bar{\partial}\tau)^{n-1} < \infty.
$$
To place this in context, recall the characterization of
${\mathbf{C}}^{n}$ by  a ``parabolic exhaustion'' $\tau$ in 
\cite{STOLL80}, \cite{DB82} where a logarithm replaces the square root in (\ref{MAESR}).
In \cite{DB82} the second named author conjectured 
that manifolds with special 
exhaustion functions  satisfying a Monge-Amp\`{e}re equation were 
affine algebraic. 
Related to this conjecture is the   
characterization of affine algebraic manifolds by exhaustions with certain 
growth conditions given by J.-P. Demailly in his 
monograph \cite{DEMAILLY}.
In fact, our task in  proving Theorem 1
 will be to show that a pair of 
 real valued functions with the growth properties 
 as in Demailly's theorem exists in any entire Grauert tube.
The holomorphic functions on ${X}$ mentioned above will have 
finite ``degree'' in the sense of Demailly's result. 

To define that pair of real 
valued functions and to find bounds for their growth 
 we make use of the 
foliation associated to the function $\tau$, the ``Monge-Amp\`{e}re 
foliation''. Properties of this  foliation  
for parabolic $\tau$ were treated in \cite{DB82} and by P.-M. Wong in 
\cite{PMW82}, but unlike 
that case, where the  singular set $\{\tau=0\}$ of the foliation is just a point, in our present 
situation the singularity set ${M}:=\{\tau=0\}$ is a totally real and 
real analytic $n$-dimensional 
submanifold.
 In light of the examples by G. Patrizio and P.-M. Wong in
 \cite{GPPMW},
 this situation was further analyzed by L. Lempert and R. Sz\H{o}ke in 
\cite{LESZ} and independently by V. Guillemin and M. Stenzel in 
 \cite{VGMS}.  
Both in \cite{LESZ} and  \cite{VGMS} it is shown that a Grauert tube 
has as a model the (co)tangent 
 bundle of ${M}$ endowed with a certain 
 ``adapted complex structure'' canonically induced from the Riemannian 
 metric that is inherited by ${M}$ from the K\"{a}hler metric 
 on ${X}$ 
 with potential $\tau$. 
 However some aspects of the foliation 
 emphasized in \cite{LESZ} 
will be the most  useful for our purposes, for instance, 
the ``holomorphic Jacobi fields'', that give outside a discrete set
a trivialization of the holomorphic tangent 
bundle along a leaf.
One of  the functions to be estimated in our proof  
 is  a Ricci curvature and thus 
we  are led to work with holomorphic $(n,0)$ forms on ${X}$. 
Using real analyticity of the foliation off the singular set  
 together with the Bott 
connection we show that
the Riemannian volume form on ${M}$ extends
to a holomorphic ``volume form'' 
 giving  a global trivialization of the canonical bundle of ${X}$. 
In the process we introduce certain ``holomorphic Fermi fields'' along 
the holomorphic parametrizations of the leaves of the foliation.
Then, the  estimates for the growth of that Ricci curvature
will essentially follow from estimating the growth of the pairing of the 
Fermi and the holomorphic Jacobi fields, leafwise, while
the compactness of the level sets of $\tau$ show that such bounds are uniform.

It follows directly from the proof above that the action of the
isometry group of the Riemannian manifold $M$ complexifies naturally
to an action on $X$ by {\em algebraic transformations}. See
\cite{RSZOKE95}. The complex manifold
$X$ also admits an anti-holomorphic involution whose fixed points are
the submanifold $M$. The algebraic embedding of Theorem 1 can be taken
to map this involution onto complex conjugation in
$\mathbf{C}^{N}$, and $M$ appears as the real points of the algebraic
variety $X$. Similar to the argument outlined above, one shows
that the metric tensor on $M$ has 
a holomorphic extension to all of $X$, and in section \ref{RATHMET}, 
Theorem 
\ref{grational}, we show that
this holomorphic form is algebraic, as well. Thus, the metric which
gives rise to the entire Grauert tube must necessarily be algebraic.

There are interesting open questions concerning existence and uniqueness of 
entire Grauert tubes, and we  will comment on them in the last section. 
For now, we state the 
following result which we prove in this paper
as an application of Theorem 1 and results of L. Lempert and R. 
Sz\H{o}ke \cite{LESZ}. 
%%%
\begin{thm}
Let ${X}$ be an unbounded Grauert tube of complex dimension
2. Then, denoting by $\simeq$ biholomorphic equivalence, either
$$
{X} \simeq 
\begin{cases}
{\mathbf C}^{*}\times{\mathbf C}^{*}
\;\; 
( {\textit{ or }}( {\mathbf C}^{*}\times
{\mathbf C}^{*} )\big{/}{\mathbf Z_{2}})
\\
{\cal  Q} \;\; ( {\textit{ or }} {\cal Q}\big{/}{\mathbf Z_2}  ) 
\end{cases}.
$$
\end{thm}
\noindent Here  ${\mathbf C}^{*}$ is the punctured complex plane
and ${\cal Q}$ is the affine $2$-quadric described above.

A version of this theorem is due independently to B. Totaro
\cite{BTotaro2}. We thank Professor Totaro for making his preprint
available to us.

\section{{\bf Holomorphic extension of the Riemannian 
volume form}}\label{EXTVOL}
We now recall some basic facts necessary to get started 
(cf. \cite{LESZ}, \cite{VGMS}, \cite{DB96}).
Let 
$({X}, \tau)$ 
be a Grauert tube of complex dimension $n$.  The Monge-Amp\`{e}re equation 
implies that $d\tau\equiv 0$ precisely at ${M}$, the zero set 
of $\tau$. In fact ${M}$ is a totally real manifold of 
dimension $n$, the fixed-point set of an anti-holomorphic involution
$\sigma$,  and ${X}$ is diffeomorphic to the tangent bundle 
of ${M }$.  
The tube has a canonical foliation, singular 
along ${M}$, the  leaves of which are Riemann 
surfaces whose tangential directions correspond to the kernel of 
$\partial \Bar{\partial}\sqrt{\tau}$. The K\"{a}hler metric on ${\cal 
X}$ with fundamental form $\sqrt{-1}\partial\Bar{\partial}\tau$
provides ${M}$
with a Riemannian metric $g$. 
By the regularity theorem of L. Lempert in \cite{LEMP}
$\tau$ is real analytic
and so is $g$ . Now, one can reconstruct the tube as the tangent 
bundle of ${M}$ endowed  with a complex structure $\mathbf{J}$ 
canonically induced from the Riemannian metric $g$ (and hence 
``adapted'' to $g$). In that model, each leaf intersects ${M}$ 
along a geodesic of $g$.

The way in which 
the Riemannian geometry of $M$
manifests itself via holomorphic objects is crucial in the  
understanding of a Grauert tube.
The following is at the heart of our argument.
%%%%%%%%%%%%%%%%%
\begin{Prop}\label{EXTENSION}
The $(n,0)$ part of the Riemannian volume form of $({M}, g)$ can be extended as a 
non-vanishing holomorphic form 
${\cal V}$ to all of ${X}$.
\end{Prop}
%%%%%%%%%%%%%%%
\begin{pf}
In local real analytic  coordinates $x_1, \cdots, x_n$ valid on a ball 
${\cal U}\subset {M}$
the Riemannian volume form has the 
expression $\sqrt{\det \tilde{g}}\, dx_1\wedge \cdots \wedge dx_n$ where 
the determinant of the matrix 
$\tilde{g}_{ij}=g(\frac{\partial}{\partial 
x_i}, \frac{\partial}{\partial 
x_j})$ is non-vanishing and real analytic. 
Thus, 
$\big{(}\sqrt{\det \tilde{g}}\, dx_1\wedge \cdots \wedge dx_n 
\big{)}^{n,0}$  
has a non-vanishing holomorphic extension
to a certain neighborhood in ${X}$ containing ${\cal U}$.
By compactness of ${M}$
 and uniqueness of analytic continuation, the volume form on ${M}$
can be  extended as a non-vanishing holomorphic form ${\cal V}$
beyond the level set $\{\tau=\epsilon \}$ for $\epsilon >0$ small enough.

We are first going to extend this holomorphic volume form 
as a non-vanishing form leaf by leaf 
to all of ${X}$ and then use real analyticity of the foliation to 
show that the result is a holomorphic form.

Fix a leaf ${\cal L}_\gamma$, intersecting  ${M}$ along a geodesic 
$\gamma$.  
Choose  a 
parametrization of $\gamma$ by arc length and  let
\begin{equation}\label{para0}
P_\gamma (x+\sqrt{-1} y) = y \dot{\gamma}(x)
\end{equation} be the  
holomorphic map 
parametrizing  ${\cal L}_\gamma$. 

\subsection{Holomorphic Fermi Fields}
We will define a special holomorphic frame along $P_\gamma$
trivializing the holomorphic vector bundle 
$P_\gamma^*(T^{1,0}({X}))$.
To do this,  
we need to recall the construction of the so-called Fermi 
coordinates.
Consider a parallel orthonormal frame along 
$\gamma$ $\{E_1(s),\cdots ,E_n(s)\}$    
    so that $E_n(s)=\dot{\gamma}(s)$, for $s$ in a neighborhood of 
    $s^0$ in ${\mathbf{R}}$. One then   
  defines the Fermi 
    coordinates $s_i$, whose domain is some  
    neighborhood ${\cal U}$ of $\gamma(s^0)$ in ${M}$, by the formulas
    \begin{equation}\label{Fcoor}
    s_i \big{(} \exp^N_{\gamma(s)} 
    (\sum_{k}^{n-1} a_k E_k(s))\big{)}
    =\begin{cases}
	a_i \textit{ for }  1 \leq i \leq n-1 \\
	s \textit{ for }  i=n \;\;,
	\end{cases}
	\end{equation}
   where $\exp^N$ is the exponential map restricted to the normal 
   bundle of $\gamma$.

Let ${z_i}$, $1\leq i \leq n$,  be the 
holomorphic extension 
of the real analytic function $s_i$ 
to a neighborhood ${\cal W} $ in ${X}$ containing 
${\cal U}$ and let
$$
Z_i := 
\frac{\partial}{\partial {z}_i},\quad 1\leq i
\leq n 
,$$ be 
the corresponding holomorphic vector fields.
The expression
$$F_i:=Z_i\circ P_\gamma $$
defines a holomorphic frame trivializing $P_\gamma^*(T^{1,0}({\cal
X}))\mid _{P_\gamma^{-1}{\cal W}}$ and by iterated use of the map
$P_\gamma$ we get  a holomorphic trivialization
of $P_\gamma^*(T^{1,0}({X}))\mid _{\cal N}$ on a neighborhood ${\cal 
N}$ of the real axis in
${\mathbf C}$.  (Here $P_\gamma$ ``unwinds'' the geodesic so that the
parallel translation of a frame $E_i$, and hence the $Z_i$'s,
correspond to well-defined sections.)

The next step is to show that the $F_i$'s can be extended from 
${\cal N}$ to all of ${\mathbf{C}}$ and that their extensions 
trivialize 
$P_\gamma^*(T^{1,0}({X}))$. Note that $F_{n}$ does extend as a non-vanishing 
holomorphic section, since from (\ref{para0}), 
(\ref{Fcoor}) and $E_{n(s)}=\dot{\gamma}(s)$
we have
the identification, for all $\omega$ in ${\cal N}$,
\begin{equation}\label{Zn}
F_{n}\circ P_{\gamma} =\frac{\partial}{\partial \omega}, 
\end{equation}
thus defining $F_{n}$ throughout ${\mathbf{C}}$.

On $\mathbf{C}$ we have the short exact sequence of vector bundles $$0
\rightarrow T^{(1,0)}(\mathbf{C}) \rightarrow P_\gamma^*(T^{1,0}(X))
\rightarrow N \rightarrow 0,$$ where $N$ is the normal bundle of
$\cal{L}_{\gamma}$ pulled back to
$\mathbf{C}$. $T^{(1,0)}(\mathbf{C})$ is spanned by $Z_n =
\frac{\partial}{\partial \omega}$, and has a connection which makes
$Z_n$ covariant constant. Similarly, $N$ has a holomorphic connection,
the Bott connection, defined locally, for $Z$ a holomorphic section of
$N$, by $$\nabla_{Z_n} Z = [Z_n, \tilde{Z}] \, {\mbox{mod}} \,
T^{(1,0)}(\mathbf{C}),$$ where $\tilde{Z}$ is a lift of $Z$ to a
holomorphic vector field on $X$ in a neighborhood. It is well known
that this is independent of the lift $\tilde{Z}$. We want to glue
these two connections together to get a holomorphic connection on
$P_\gamma^*(T^{1,0}(X))$. To do this, it suffices to show that the
short exact sequence splits holomorphically in a natural way. Indeed,
off the real axis, the one form $\partial u$ is holomorphic and its
vanishing defines a complementary bundle $E := \, {\mbox{Ker}} \; \partial
u$ to $T^{(1,0)}(\mathbf{C})$, which enables us to define a connection
on $P_\gamma^*(T^{1,0}(X))$ as the sum of the two pieces already
discussed on the two summands (lifting the Bott connection to the
complement of $T^{(1,0)}(\mathbf{C})$). We need to show that this
splitting of the bundle $P_\gamma^*(T^{1,0}(X))$ can be continued
holomorphically across the real axis.

To do this, we note that since $u = \sqrt{\tau}$, we have that
$\partial u = \frac{1}{2u} \partial \tau$, away from $M$. On
the upper half plane in $\mathbf{C}$, we  have that $u = y_n = Im \;
\omega$; on the lower half plane, $u = - y_n = - Im \; \omega.$ Let us
check that each of the functions $-2 \sqrt{-1} <\partial u, Z_j>$ on
the upper half plane near the real axis extends holomorphically across
the real axis. More precisely, we wish to show they extend
holomorphically past the real axis and $$-2 \sqrt{-1} <\partial u,
Z_j> = \delta_{j,n}$$ along the axis, where $\delta_{j,n}$ is the
Kronecker $\delta$. Note first that for $j = n$, the conclusion is
obvious, since $<\partial u, Z_n> \equiv -\frac{\sqrt{-1}}{2}.$ For
the others, we first note that $$ <\partial u, Z_j> = \frac{1}{2u}
Z_j(\tau) = \frac{1}{2y_n} \frac{\partial
\tau}{\partial z_j}.$$ Since $\tau$ has a minimum along $M$, $\frac{\partial
\tau}{\partial z_j} \equiv 0$ along $\gamma$ (or along the real axis
in $\mathbf{C}$). Restricted to $\cal{L}_{\gamma}$ (pulled back to
$\mathbf{C}$) where $y_n = 0$ is a real analytic defining equation for
$\gamma$ (the real axis), we conclude that $y_n$ divides $\frac{\partial
\tau}{\partial z_j}$ real analytically in a neighborhood of
$\gamma$. In particular, the extended function has
value along $\gamma$ given by $\frac{1}{2} \frac{\partial^2
\tau}{\partial y_n \partial z_j}$. But since $\frac{\partial
\tau}{\partial z_j} \equiv 0$ along $\gamma$, we have
$\frac{\partial^2 \tau}{\partial x_n \partial z_j} \equiv 0$ along
$\gamma$, and hence $$\frac{\partial^2 \tau}{\partial y_n \partial
z_j} = -2\sqrt{-1} \frac{\partial^2 \tau}{\partial \bar{z}_n \partial
z_j},$$ along $\gamma$. Finally, $\frac{\partial^2 \tau}{\partial
\bar{z}_n \partial z_{j}} \equiv 0$ along $\gamma$ because the 
$\tau$-metric restrictsto give our original metric on $M$ and because 
of the orthogonality properties of the Fermi fields along $\gamma$.

Note that by Schwarz reflection, the extension of $-2\sqrt{-1}
\partial u$ to the lower half plane agrees with $2\sqrt{-1} \partial
u$ there. Thus, we have extended the holomorphic complement to
$T^{(1,0)}(\mathbf{C})$ holomorphically across the real axis, and
therefore constructed a holomorphic connection on the bundle
$P_\gamma^*(T^{1,0}(X))$ over all of $\mathbf{C}$.

On $ {\cal N} \cap \gamma$, the $F_{i}$'s lie in the subbundle $E$
complementary to $T^{(1,0)}(\mathbf{C})$. Since they are holomorphic
sections, they must lie in $E$ on all of $\cal{N}$. Then the $F_i$'s
can be taken to be the lifts to $P_\gamma^*(T^{1,0}(X))$ of the
sections of $N$ on $\cal{N}$ given by the Fermi coordinate fields
$Z_1, \cdots, Z_{n-1}$, and these $F_i$ are covariant constant since
they are coordinate fields. Finally, this allows us to extend them so
that together with $F_{n}$ they form a holomorphic frame trivializing
$P_\gamma^*(T^{1,0}({X}))$.

To complete the proof of Proposition \ref{EXTENSION}, let us consider
the dual forms of the $F_{i}$, i.e., the sections $F_{i}^{*}$ of the
holomorphic bundle $P_{\gamma}^{*}(T^{1,0}({X}))^{*}$. Since along
${\cal N}$
$$F_{1}^{*}\wedge \cdots \wedge F_{n}^{*}= {\cal V}\circ P_{\gamma},$$
we can now extend ${\cal V}$ from ${\cal W}\cap {\cal L}_{\gamma}$ to 
all of ${\cal L}_{\gamma}$. 

Because ${X}\setminus {M}$ is foliated by the ${\cal 
L}_{\gamma}$'s, 
the construction above performed for each geodesic $\gamma$  
defines a form  
${\cal V}$ for every point of ${X}$. In principle we 
only  
know this form to be holomorphic in the region 
$\{0\leq \tau < \epsilon \}$, 
but because of real analyticity of $\tau$ we have the following   

Claim: The leaf-wise extension ${\cal V}$ is holomorphic on 
${X}$.

To see this, it suffices to observe that the $(n,0)$-form $\cal{V}$ is
real analytic on $X$, since we already know that it is holomorphic in
a neighborhood of $M$, and would be holomorphic on all of $X$ by
analytic continuation. On the other hand, the geodesic $\gamma$ (and
consequently the leaf $\cal{L}_{\gamma}$) depends real analytically on
$2n-2$ real parameters $s = (s_1, \cdots, s_{2n-2})$, as do the
associated Fermi coordinates $x = x(s)$ and their complexifications $z
= z(s)$, and the corresponding Fermi fields $Z_i(s)$,
locally. Finally, since the Monge-Amp\`ere solution $u$ is
$\cal{C}^{\omega}$ in all arguments off $M$, the complementary bundles
$E(s)$ also vary real analytically. The extensions of the $Z_i(s)$
along the leaves are achieved by solving systems of holomorphic
ordinary differential equations which depend real analytically on the
parameters $s$, and with real analytically varying initial
conditions, hence they depend real analytically on all
parameters. This concludes the proof of Proposition \ref{EXTENSION}.

\end{pf}

%%%%%%%%%%%%%%%%%%%%%%%%%%%%%%%%%%%%%%%%%%%%%%%%%%%%%%%%%%%%%%%%%%%%%%%%%%%%%%%%%%%%%%%%%%%%%
\section{{\bf J.-P. Demailly's characterization of affine 
varieties}}
Let's recall the following very pretty result, 
Theorem $9.1^\prime$ in \cite{DEMAILLY}, 
stated here in less 
generality than the original version, just to the 
degree needed for  our purposes.
%%%%%%%%%%%%%%%%%%%%%%%%%%%%%%%%%%%%%

{\bf Theorem:} (Demailly \cite{DEMAILLY})
{\em Let ${X}$ a complex manifold of complex dimension $n$ with
finite dimensional  cohomology groups. Assume that there is a  
smooth strictly plurisubharmonic exhaustion 
function $\phi$ on ${X}$ with  
$$\int_{X} (\sqrt{-1}\partial\Bar{\partial}\phi)^n <\infty.
$$
Let $R(\beta)$ be the Ricci form of the K\"{a}hler metric
induced by the potential $e^\phi$ and put 
$\beta=\sqrt{-1}\partial \Bar{\partial}e^{\phi}$.
If there are non-negative constants $c$ and $c^{\prime}$,  and
 a $C^{2}$ function $\psi$ 
 so that
$$
R(\beta)+\sqrt{-1}\partial \Bar{\partial} \psi \geq 0, \;\;
{\textit{ and }}
\;\;
\int_{X}e^{c\psi-c^{\prime} \phi}\beta^{n}<\infty$$
then ${X}$ is biholomorphic to an affine algebraic variety.}

The algebraic embedding of ${X}$ into ${\mathbf{C}}^{N}$, $N>>1$,  
is via holomorphic functions in the algebra 
$A_{\phi}$ consisting of  ``$\phi$-polynomials'',  
holomorphic functions $f$ in ${X}$ such that
their ``degree''$$ 
\lim \sup_{r\to \infty}
\frac{1}{r}
\int_{\phi=r}
\log _{+}\mid f\mid  d\phi \wedge (\partial \Bar{\partial}\phi)^{n-1}
 < \infty.
$$
In this section we propose for a 
given Grauert tube $\big{(}{X}, \tau\big{)}$ 
of infinite radius functions $\phi$
and $\psi$ that will 
meet the conditions in Demailly's theorem. 
The proof that $\phi$ indeed satisfies the required growth conditions 
occupies the next two sections. 

The $\phi$ will be a function solely of $\tau$ and 
the $\psi$ will be defined in terms of the holomorphic extension of 
the Riemannian volume form from Proposition \ref{EXTENSION}.
The choices made are not unique but are in fact 
quite natural (cf. Remark \ref{choices}).
%%%%%%%%%%%%%%%%%%%%%%%%%%%%%%%%%
\begin{Prop}
Let ${X}$ be an unbounded tube with Monge-Amp\`{e}re exhaustion 
function $\tau$. Then 
\begin{equation}\label{defphi}
    \phi=\ln (1+\cosh \sqrt{\tau})
\end{equation}
    is strictly plurisubharmonic
and $\int_{X} (\sqrt{-1}\partial \bar{\partial} \phi)^{n} <\infty.$ 
\end{Prop}
%%%%%%%%%%%%%%%%%%%%%%%%%%%
\begin{pf}
We have $\phi=f(\tau)$ and $\partial \bar{\partial} \phi 
= 
f^{\prime \prime}(\tau) \partial \tau \wedge \bar{\partial}\tau 
+
f^{\prime} (\tau) 
\partial\bar{\partial}\tau
$
where
$$
f^{\prime}(\tau)=\frac{1}{1+\cosh \sqrt{\tau}}
\frac{ \sinh \sqrt{\tau} }{ 2 \sqrt{\tau}}, \hspace*{2mm}
f^{\prime \prime}(\tau) =
\frac{1}{4\sqrt{\tau} (1+\cosh \sqrt{\tau})}
\Big{(}1-\frac{\sinh \sqrt{\tau}}{\sqrt{\tau}} 
\Big{)}.$$
Strict plurisubharmonicity easily follows
from the inequalities $    f^{\prime}(\tau) >0$ and  
$2 \tau f^{\prime \prime}(\tau) +f^{\prime}(\tau)>0$ and 
the splitting discussed in the proof of Proposition \ref{EXTENSION},
$T^{1,0}({X})|_{{X}\setminus {M}}
={\mathbf{C}} \xi \oplus 
{\textit{ Kernel }} \partial \tau$, with $\xi$ in the tangential direction
to each leaf.

On the other hand,
from the Monge-Amp\`{e}re equation we get
\begin{equation}\label{maddb}\partial \tau \wedge \bar{\partial}\tau
\wedge 
 (\partial\bar{\partial}\tau )^{n-1}
=\frac{2\tau}{n}
 (\partial\bar{\partial}\tau )^{n}
\end{equation}
and by letting
$$\delta_{0} (\tau) := \frac{1}{2^{n}}\Big{(}\frac{\sinh 
\sqrt{\tau}}{\sqrt{\tau}}\Big{)}^{{n-1}},
$$
we have 
$$
(\partial\Bar{\partial}\phi)^n=
\frac{\delta_{0}(\tau)}{(1+\cosh \sqrt{\tau})^{n}}
(\partial \Bar{\partial} \tau)^n.
$$
Putting  ${X}^\prime :=\{1\leq \tau < \infty \}$   
it follows that there is a $C>0$ so that
\begin{equation}\label{intfin}  
  \int_{{X}}
  (\sqrt{-1} \partial \Bar{\partial} \phi)^n
  \le C
 \int_{{X}^\prime}
 e^{-\sqrt{\tau}}
(\sqrt{-1} \partial\Bar{\partial}\tau)^n.
\end{equation}
But under the identification of ${X}$ with the tangent bundle
of ${M}=\{ \tau = 0 \}$ endowed with the Riemannian metric $g$,
 $\tau$ corresponds to the $g$-length squared in $T{\cal 
M}$ and $\sqrt{-1}\partial \Bar{\partial}\tau$ to the symplectic 
$2$-form induced via $g$ by the tautological $2$-form in the 
cotangent bundle $T^*{M}$. Thus, we have 
\begin{equation*}
\begin{align}
\int_{{X}^\prime}
 e^{-\sqrt{\tau}}
  (\sqrt{-1} \partial\Bar{\partial}\tau)^n = &
 \,n\int_{{X}^\prime}
 \frac{e^{-\sqrt{\tau}}}{2\tau}
 \partial{\tau}\wedge \Bar{\partial}\tau \wedge 
(\sqrt{-1} \partial\Bar{\partial}\tau)^{n-1}\notag \\
= & C^\prime {\textit{ Volume  of }} ({M},g)
\int_1^\infty 
\lambda^{n-2} e^{-\sqrt{\lambda}} 
d\lambda \le \infty \notag 
\end{align}
\end{equation*}
for a constant $C^\prime$ depending on $n$,
showing that 
$
  \int_{X} (\sqrt{-1}\partial \Bar{\partial} \Phi)^n$
is finite.
\end{pf}
%%%%%%%%%%%%%%%%%%%%%%%%%%%%%%%%%%%%%%%%%%%%%%%%%%%%%%%%%%%%%%%%%%
Next, in order to find a candidate for $\psi$,  we  consider
the Ricci form $R(\beta)$ of the K\"{a}hler metric $\beta$
induced by the potential
\begin{equation}
    e^{\phi}=1+\cosh \sqrt{\tau}.
\end{equation}
We have,
\begin{equation}\label{betametric}
\partial \Bar{\partial}(1+ \cosh \sqrt{\tau})=
\frac{\sinh \sqrt{\tau}}{2\sqrt{\tau}}
\Big{(}
\frac{1}{2\tau}
(\frac{\sqrt{\tau}}{\tanh\sqrt{\tau}}
- 1)
\partial \tau \wedge\Bar{\partial}\tau 
+\partial\Bar{\partial} \tau \Big{)},
\end{equation} 
and hence, using 
(\ref{maddb}),
\begin{equation}\label{volcosh}
(\partial \Bar{\partial} (1+ \cosh \sqrt{\tau}))^n=
 \delta_{0}\, \cosh \sqrt{\tau}\; 
 (\partial\Bar{\partial}\tau)^n.
\end{equation}
%%
%%%
Let ${\cal V}$ be the non-vanishing holomorphic
$(n,0)$-form on ${X}$ from Proposition \ref{EXTENSION} 
extending the Riemannian volume 
form on ${M}$, and let 
$\delta$ be the function 
$$\delta\colon {X}\to {\mathbf R}>0$$
globally defined by
\begin{equation}\label{gvolddrho} 
{\cal V} \wedge \Bar{{\cal V}}=
\frac{(-1)^{\frac{n^2}{2}}}{n!}
\delta  (\sqrt{-1}\partial \Bar{\partial} \tau)^{n}.
\end{equation}
By (\ref{volcosh}) and (\ref{gvolddrho}), the Ricci form $R(\beta)$ 
has  globally the expression
$$R(\beta)=\frac{1}{\sqrt{-1}}
\partial\Bar{\partial}
\ln \big{(}\frac{\delta_{0}}{\delta}\,
\cosh \sqrt{\tau} \big{)},
$$
and 
if we simply let
\begin{equation}\label{defpsi}
\psi:= \ln \big{(}\frac{\delta_{0}}{\delta}\,
\cosh \sqrt{\tau} \big{)}, 
\end{equation}
then
\begin{equation}\label{Dem1}
    R(\beta)+\sqrt{-1}\partial\Bar{\partial}\psi 
    =0. 
\end{equation}
Consider the following
%%
%%%%%%%%%%%%%%%%%%%%%%%%%%%%%%%%%%%%%%%
\begin{Prop}\label{ESTIMATEM}
 $ \exists$ constants  $C>0$ and $k>0$ so that
\begin{equation}\label{ESTIMATEMf}
1/\delta \le C e^{k\sqrt{\tau}}.
\end{equation}
\end{Prop}
%%%%%%%%%%%%%%%%%%%%%%%%%%%%%%%%%%%%%%%
Since clearly we can find a constant $c>0$ so that 
$
\cosh \sqrt{\tau}\; \delta_0 (\tau)< ce^{(n+1)\sqrt{\tau}},
$ 
then
Proposition \ref{ESTIMATEM} 
has the following (cf. (\ref{defphi}), (\ref{defpsi}))
\begin{Cor}
 $ \exists$  constants $A>0$ and $B>0$ so that
\begin{equation}\label{Dem2}
    \psi \leq A\phi +B.
    \end{equation}
\end{Cor}
\subsubsection{Proof of Theorem 1}
\begin{pf} We prove Proposition \ref{ESTIMATEM}, and hence (\ref{Dem2}),
 in the next two 
sections. Assuming those results, 
by (\ref{Dem1}) and (\ref{Dem2}) the functions $\phi$ and 
$\psi$  meet 
the conditions of Demailly's theorem 9.1' in \cite{DEMAILLY}. 

Indeed, given 
(\ref{Dem2}) it is immediate to check that
there is a constant $c^{\prime \prime}>0$ so that 
$\int_{X}e^{\psi - c^{\prime}\phi}\beta^{n}$ is bounded above by
$$
c^{\prime \prime }\int_{X}
e^{(n+1)\sqrt{\tau}} \big{(} 1+\cosh \sqrt{\tau}\big{)}^{A-c^{\prime}} 
(\sqrt{-1}\partial \Bar{\partial}\tau)^{n},$$
a finite number for $c^{\prime}>>A$. 

 The  finiteness 
of the cohomology groups is satisfied since ${X}$ is 
diffeomorphic to the tangent bundle of  
the  compact manifold $\{\tau =0\}$.
\end{pf}

\begin{Rem}\label{choices}
The choices of $\phi$ and $\psi$ made above, although not unique, are natural 
($\phi$ was already noticed to have the ``finite volume'' property in 
\cite{GPPMW}). 
\end{Rem}

For example, consider the affine quadric ${\cal
Q}_{n}=\{\sum_{i=1}^{n+1}\omega^{2}_{i}=1\}
\subset {\mathbf C}^{n+1}\subset {\mathbf CP}^{n+1}$.
By \cite{GPPMW}, 
$\sqrt{\tau}=\cosh ^{-1} \sum \mid \omega_{i} \mid ^{2}$ 
and hence
$\phi$ is the restriction of the potential for the 
Study-Fubini metric in ${\mathbf CP}^{n+1}$.
A computation will show that $\delta=\delta_{0}$ and 
$\psi=\ln\sum \mid \omega_{i} \mid ^{2}$.  
Note that in general $\delta$ is not a function 
of $\tau$ alone. It is so precisely when the Riemannian metric induced 
in $\{\tau=0\}$ is {\em harmonic} (cf. \cite{RMA2}).

We devote the next two  sections 
to the  proof of 
Proposition \ref{ESTIMATEM}, which will complete the proof of Theorem 1.
%%%%%%%%%%%%%%%%%%%%%%%%%%%%%%%%%%%%%%%%%%%%%%%%%%%%%%%%%%%%%%%%%%%%%%%%%%%%
\section{{\bf Proof of Proposition \ref{ESTIMATEM}:
                   Part 1,  $\delta $ along a leaf.}}\label{PFP1}
%%%%%%%%%%%%%%%%%%%%%%%%%%%%%%%%%%%%%%%%%%%%%%%%%%%%%%%%%%%%%%%%%%%%%%%%%%%
In this section once again we identify as complex manifolds 
the Grauert tube ${X}$ with 
$(T{M}, \mathbf{J})$, the 
tangent bundle of 
${M}=\{\tau =0\}$ 
endowed with the adapted complex structure ${\mathbf{J}}$.
Each geodesic $\gamma$ of $({M}, g)$
 defines a leaf ${\cal L}_\gamma$, 
 which, for a chosen  parametrization of $\gamma$ by arc length,
is given a holomorphic parametrization 
$P_{\gamma}\colon {\mathbf C} \to{M}$
defined  by 
\begin{equation}\label{pgamma}
    P_\gamma (\omega) = P_\gamma ( x+\sqrt{-1}y) =y \, \dot{\gamma}(x).
\end{equation}
Our goal is to show that there are constants $C>0$ and $k>0$ such that
$1/\delta <Ce^{k\sqrt{\tau}}$. We will accomplish this in the next
section. In the present section we first obtain a convenient
expression of $(1/\delta)\circ P_\gamma$ for each parametrized
geodesic $\gamma$, given in Proposition \ref{PDELTA}.  (In fact we will
later be able to show the inequality (\ref{strineq}), which is
stronger than needed here.)

Choose an orthonormal frame of $T_{\gamma(0)}{M}$ and extend 
it along $\gamma$ by parallel translation to obtain
the frame $
E_1{(s)},\cdots E_{n}(s)$,
which at each $s$ is  an orthonormal basis of $T_{\gamma(s)}{M}$.
Denote by the same symbols the vector fields on
the real axis along $P_{\gamma}$.

Consider now the corresponding 
holomorphic Fermi fields along $P_{\gamma}$ as we defined them in 
the proof of Proposition \ref{EXTENSION}, 
\begin{equation}\label{Fframe}
F_{1}(\omega),\cdots F_{n}(\omega),
\end{equation}
and recall, for later reference,  that for all $-\infty < x <\infty$, $1\leq 
i \leq n$, 
\begin{equation}\label{FFRA}
F_{i}(x)=E^{1,0}_{i}(x).
\end{equation}
For all $\omega \in {\mathbf C}$, again from our construction 
in the proof of Proposition \ref{EXTENSION}, 
$$
{\cal V}=\Pi_{i=1}^{n}F_{i}^{*}(\omega),
$$
where the $F^{*}_{i}$'s are the forms dual  to the $F_{i}$'s
and where, to simplify notation, 
 we also denote by ${\cal V}$ the holomorphic section of 
$P_{\gamma}^{*}\big{(}(T^{n,0}({X})^{*}\big{)}$ 
determined by the holomorphic 
volume form ${\cal V}$.
Thus
$$
\Pi_{i=1}^{n} F^{*}_{i} \wedge \Bar{F}^{*}_{i}=
(-1)^{\frac{n(n-1)}{2}}
 {\cal V} \wedge \Bar{\cal V},
$$
i.e., 
\begin{equation}\label{volumefermi}
\Pi_{i=1}^{n} F^{*}_{i} \wedge \Bar{F}^{*}_{i}=
\frac{\delta}{n!} (\partial\Bar{\partial}\tau)^{n}.
\end{equation}
To take advantage of 
equation (\ref{volumefermi})
we will resort to two more sets of vector fields holomorphic along 
$P_{\gamma}$. The basic reference here is L. Lempert and R. 
Sz\H{o}ke's work \cite{LESZ}, where they introduced these vector fields 
to describe the adapted complex structure.

Let $\xi^{0}_{i}(s)$ and 
$\eta^{0}_{i}(s)$, $1\leq i \leq n$, be the 
Jacobi fields along $\gamma$ satisfying ``Lagrangian'' 
initial conditions such as
\begin{equation}\label{initialc}
\eta_{i}^{0}(0) =\nabla_{\dot{\gamma}(0)}\xi^{0}
_{i}(0)= 0,\quad 
\xi^{0}_{i}(0) =\nabla_{\dot{\gamma}(0)}\eta^{0}_{i}(0)= E_i (0), 
\end{equation}
where $\nabla$ is the Levi-Civita  connection on $({M}, g)$.
Recall that these vector fields are solutions of the differential equation
\begin{equation}\label{Jacobieq}
\nabla^{2}_{\dot{\gamma}(x)}\xi^{0}_{i}(x)+
 R_{\gamma}(x)\xi^{0}_{i}(x)=0,
\end{equation}
(and similarly for the $\eta_{i}$'s), where 
$R_{\gamma}(x)\colon T_{\gamma(x)}{M}\to T_{\gamma(x)}{M}$
is the curvature operator along $\gamma$ given by 
$v\mapsto R(v, \dot{\gamma}(0))\dot{\gamma}$ with $R$ the curvature 
tensor of $({M}, g)$. They 
determine along $P_{\gamma}$ vector fields $\eta_{i}$ and $\xi_{i}$
that are invariant by the geodesic flow and fiberwise rescaling.
It is shown in \cite{LESZ} that the $(1,0)$ parts of these invariant 
fields
\begin{equation}\label{hjvf}
\xi_{i}^{1,0}=\frac{1}{2}(\xi_{i}-\sqrt{-1}{\mathbf{J}}\xi_{i}), 
\qquad 
\eta_{i}^{1,0}=\frac{1}{2}(\eta_{i}-\sqrt{-1}{\mathbf{J}}\eta_{i})
\end{equation}
are holomorphic sections of the holomorphic bundle over 
${\mathbf{C}}$ $P_{\gamma}^{*} \big{(}T^{1,0}({X})\big{)}$ (these are 
the ``holomorphic Jacobi 
fields'' referred to above).
Again using the same symbols $\eta^{0}_{i}$ and $\xi_{i}^{0}$ 
to denote
the induced fields on  the real axis along $P_{\gamma}$ we have
$$(\eta^{0})^{1,0}_{i}(x)=\eta^{1,0}_{i}(x),
\qquad 
(\xi^{0})^{1,0}_{i}(x)=\xi^{1,0}_{i}(x).
$$
Now, on ${\mathbf R}\setminus S$, where 
%%%%%%
$$S:=\{x \in {\mathbf{R}} \mid  \xi_1^0(x)\wedge\cdots 
\wedge\xi_n^0(x)=0\},$$ 
%%%%
a discrete set,
for the $\xi_{i}$'s 
are solutions along the real axis of 
the second order differential equations (\ref{Jacobieq}),  
we also have
\begin{equation}\label{etahxi}
 \eta^{0}_{i}(x)=\sum_{k=1}^{n}h^{0}_{ki}(x) \xi^{0}_{k}(x).
 \end{equation}
It follows that along $P_{\gamma}$ there are meromorphic functions
$h_{ki}(\omega)$ extending $h^{0}_{ki}(x)$ to ${\mathbf C}$
so that
\begin{equation}\label{etahxidef}
\eta^{1,0}_{i}(\omega)=\sum_{k=1}^{n}h_{ki}(\omega)\xi^{1,0}_{k}(\omega),
\end{equation}
which implies
\begin{equation}\label{jdef}
\eta_{i}(\omega) =\sum_{k=1}^{n}
\Re h_{ki}(\omega) \; \xi_{k}(\omega) +\sum_{k=1}^{n}
\Im h_{ki}(\omega) \;{\mathbf{J}}\xi_{k}(\omega).
\end{equation}
We will need to recall from \cite{LESZ} 
more facts about these meromorphic functions 
later, but for now  
we note that the volume form of the 
K\"{a}hler metric given by $\tau$ can be expressed as 
$$\frac{1}{n!}(\partial\Bar{\partial}\tau)^{n}=
\det M \;\Pi_{i=1}^{n}(\xi_{i}^{1,0})^{*}\wedge (\xi_{j}^{0,1})^{*}
$$
where, again ``$*$'' denotes dual forms and 
\begin{equation}\label{MIJD}
    M_{ij}(\omega) =
    \partial\Bar{\partial}\tau (\xi^{1,0}_{i}, \xi^{0,1}_{j}). 
\end{equation}
Furthermore, since by the proof of Proposition \ref{EXTENSION} 
the holomorphic Fermi fields along $P_{\gamma}$ are 
$\mathbf{C}$-linearly independent on all of ${\mathbf{C}}$ 
we have
\begin{equation}\label{AKI}
\xi_{i}^{1,0}(\omega) = \sum_{k=1}^{n} A_{ki}(\omega) F_{k}(\omega),
\end{equation}
 where the 
$A_{ki}(\omega)$ are holomorphic functions
on ${\mathbf C}$ so that on ${\mathbf{C}}\setminus S$ 
$$\det A (\omega) \neq 0 .$$ 
%%%%%%
 We thus conclude, noting that the dual forms corresponding to  the fields in 
 (\ref{AKI}) are 
 related by the matrix $(A^{-1})^{T}$,  that
 %%%%%%%%%%%%%%%%%%%%%%
\begin{equation}\label{espest}
\frac{1}{\delta} \circ P_{\gamma}= 
\frac{\det M (\omega )}{\mid \det A 
(\omega) \mid ^{2}}\quad {\textit{ for all }}\Im \omega\neq 0.  
\end{equation}
%%%%%%%%%%%%%%%%%%%%%%%
Because $\tau$ is an exhaustion function we only need to 
estimate  $1/\delta$ outside a relatively compact neighborhood 
of ${M}$, say on $\tau\geq1$.
In particular, for any geodesic $\gamma$,  it will suffice 
to bound  $(1/\delta) \circ P_{\gamma}$ for 
 $y\geq1$.
  
We will first deal with an estimate of the  numerator
of (\ref{espest}), our present goal being to 
express it in terms of the functions $h_{ki}(\omega)$, for which 
we know how to obtain bounds. 

With that in mind, we
 will rewrite the coefficients $M_{ij}(\omega)$ (\ref{MIJD}) in the form 
 (\ref{MIJ}).
 %%%%%%%%%%%%
 Note that
at each point $\omega=x+\sqrt{-1}y$ in the leaf with $y > 0$
the span of the vector fields $\{\xi_1,\cdots , \xi_n \}$
 and that of the $\{\eta_1,\cdots , \eta_n \}$ are two complementary 
 Lagrangian subspaces of $T_{P_{\gamma}\omega}\big{(}T{M}\big{)}$ with,
 \begin{equation}\label{lagrange0}
 \sqrt{-1}\partial\Bar{\partial}\tau (\xi_{i}, \xi_{j}) =
 \sqrt{-1}\partial\Bar{\partial}\tau (\eta_{i}, \eta_{j})=0
 \end{equation}
\begin{equation}\label{lagrange1}
 \sqrt{-1}\partial\Bar{\partial}\tau (\xi_{i}, \eta_{j}) =y\delta_{ij}.
 \end{equation}
Equations (\ref{lagrange0}) and (\ref{lagrange1}) are a consequence 
of the initial conditions (\ref{initialc})
and the equivariance of the $\xi_{i}$, the $\eta_{i}$'s and the 
symplectic form.
Indeed, at each point $\omega$ in the leaf, 
\begin{equation}\label{KPI}
K\xi_{i}(\omega) =\nabla_{P\omega}\xi^{0}(\pi(P_{\gamma}\omega)), \quad 
\pi_{*}\xi_{i}(\omega)= \xi^{0}_{i}(\pi(P_{\gamma}\omega)), 
\end{equation}
and similarly for the $\eta_{i}$'s, 
where $K\colon T\big{(}T{M}\big{)}\to T{M}$ is the connection map
and $\pi\colon T{M}\to {M}$ the natural projection.
Moreover,  $\sqrt{-1}\partial\Bar{\partial}\tau$ 
which has the expression, always under the tangent bundle 
identification,  
\begin{equation}\label{SYMPLECTIC}
\sqrt{-1}\partial\Bar{\partial}\tau 
= g\circ (K\otimes \pi_{*}-\pi_{*}\otimes K),
\end{equation}
is invariant by the 
geodesic flow, and homogeneous of degree one with respect to multiplication  
along the fibers of $T{M}$. Thus, equations (\ref{lagrange0}) and 
(\ref{lagrange1}), valid by (\ref{initialc}) at 
$P_{\gamma}^{-1}(\dot{\gamma}(0))$, are so throughout the whole leaf 
${\cal L}_{\gamma}$.

It follows that, 
\begin{equation*}
M_{ij}(\omega)=\frac{1}{2}\sqrt{-1}
\partial\Bar{\partial}\tau (\xi_{i}(\omega), 
{\mathbf{J}} \xi_{j}(\omega)),
\end{equation*}
and from 
(\ref{jdef}), (\ref{lagrange0}) and (\ref{lagrange1}),
 for all $\omega $ with $\Im \omega>0$, 
%%%
 \begin{equation}\label{MIJ}
M_{ij}(\omega)=\frac{1}{2} \Im \omega  (\Im h)^{{-1}}_{ij}(\omega).
\end{equation}
%%%
%%%%%%%%%%%%%%%%%%  %%
Now we focus on the denominator of (\ref{espest}), again our objective 
being to express it in terms of the functions $h_{ki}(\omega)$, which 
we do by obtaining (\ref{detAcplane}).
%%%%%%%%%%%%%%%%%  %%

Equation (\ref{AKI}) reads along the real axis, by (\ref{FFRA}), 
\begin{equation}\label{AKIRA}
(\xi_{i}^{0})^{1,0}(x)= \sum_{k=1}^{n} A_{ki}(x)E^{1,0}_{k}(x).
\end{equation}
But, since $T_{\gamma(x)}{M}$ is spanned by the $E_{k}(x)$'s while 
along the zero section $s\colon {M}\to T{M}$ 
we have $T\big{(}T{M}\big{)}|_{M}=
s_{*} (T{M})\oplus {\mathbf{J}} s_{*}( T{M}),
$
it follows that for $-\infty < x < \infty$
\begin{equation}\label{AKIRA2}
\xi_{i}^{0}(x)= \sum_{k=1}^{n} A_{ki}(x)E_{k}(x).
\end{equation}
The same arguments concerning the $\xi_{i}$'s 
leading  to (\ref{AKI}), (\ref{AKIRA}) and (\ref{AKIRA2})
 apply to the $\eta_{i}$'s so that
for $\omega$ in ${\mathbf{C}}$ 
\begin{equation}\label{BKI}
\eta_{i}^{1,0}(\omega) =\sum_{k=1}^{n}B_{ki}(\omega) F_{k}(\omega)
\end{equation}
where $B_{ki}$ are holomorphic on ${\mathbf{C}}$
and moreover along the real axis
\begin{equation}\label{BKIRA}
\eta^{0}_{i}(x) =\sum_{k=1}^{n} B_{ki}(x) E_{k}(x).
\end{equation}
Of course, from (\ref{etahxi}), (\ref{AKI}) and (\ref{BKI}), 
for all $\omega$ outside the discrete subset 
$S\subset {\mathbf{R}}$ we have
\begin{equation}\label{BAF}
B_{ki}(\omega)= \sum_{s=1}^{n}A_{ks}(\omega) h_{si}(\omega),
\end{equation}
which gives,  upon restriction to ${\mathbf{R}}\setminus S$, 
double differentiation 
and use of the Jacobi equations (\ref{Jacobieq})
together with the fact that the $E_{i}$'s are parallel along 
$\gamma$, 
the differential equation (cf. Proposition 6.11 in \cite{LESZ}), now 
in matrix notation,
\begin{equation}\label{DIFFEQ}
    2 A^{\prime}(x) h^{0\,\prime}(x)+A(x) h^{0\,\prime\prime}(x)=0.
\end{equation}
On the other hand, from (\ref{lagrange0}), (\ref{KPI}) and 
(\ref{SYMPLECTIC}) it follows that (cf. Proposition 6.10 in \cite{LESZ})
$\big{(}A(x)^{T}A(x)-A(x)A(x)^{T}\big{)}^{\prime}=0$.
Using this back in (\ref{DIFFEQ}) 
and taking into account 
the initial conditions for 
the $\xi^0_i$'s, and $\eta^0_i$'s
(\ref{initialc}) for the equations (\ref{Jacobieq}) which imply
$A(0)=h^{0\,\prime}(0)={\mathbf{I}}$, 
gives 
for all $\omega$ in ${\mathbf{C}}\setminus S$,
\begin{equation}\label{ATAHPRIME}
A^{T}(\omega)A(\omega)=A(\omega)A^{T}(\omega)=(h^{\prime}(\omega))^{-1}.
\end{equation}
In particular, 
\begin{equation}\label{detAcplane}
\big{(}\det  A(\omega) \big{)}^{-2} = \det h^{\prime}(\omega) .
\end{equation}

%%%%%%%%%%%%%%%%%%%%%%%%%
To summarize,  in virtue of (\ref{espest}), (\ref{MIJ}) and  
(\ref{detAcplane}), we have shown
%%%%%%%%%%%%%%%%%%%%%
\begin{Prop}\label{PDELTA}
    Let $\gamma$ be an arc length  parametrized geodesic, $P_\gamma$ 
    the map (\ref{pgamma}) and $h_\gamma$ the matrix with components 
    given by (\ref{etahxidef}) corresponding 
    to $\gamma$. Then 
for all $y= \Im \omega >0$,
\begin{equation}\label{PDELTAEST}
\frac{1}{\delta(P_{\gamma}\,\omega)} =  
y^n 
\frac{ \mid  \det h_\gamma^{\prime}(\omega ) \mid }
     {\det (\Im h_\gamma) (\omega) }.
\end{equation}
\end{Prop}
%%%%%%%%%%%%%%%%%

%%%%%%%%%%%%%%%%%%%%%%%%%%%%%%%%%%%%%%%%%%%%%%%%%%%%%%%%%%%%%%%%%%%%%%%%%%%%%%%
\section{{\bf Proof of Proposition \ref{ESTIMATEM}:  
                         Part 2, growth estimate for 
                         $\delta$.}}\label{PFP2}
%%%%%%%%%%%%%%%%%%%%%%%%%%%%%%%%%%%%%%%%%%%%%%%%%%%%%%%%%%%%%%%%%%%%%%%%%%%%%

In this section we will use (\ref{PDELTAEST}) to estimate the 
growth of $\delta$ in terms of $\tau$. 
This will complete the proof 
of Proposition \ref{ESTIMATEM}.

We show our bounds  along the 
same lines of R. Sz\H{o}ke's work in \cite{RSZOKE95}. 
Consider the upper half-plane 
 $$
 {\mathbf{C}}^+=\{\omega =x+\sqrt{-1}y 
 \in {\mathbf{C}}, x,y \in {\mathbf{R}} \mid  y >0 \}
 $$
 and  the space of symmetric complex valued $n\times n$ matrices
with positive semi-definite
imaginary part,
$$
 {\cal H}^n=\{M \in M^n_{\mathbf{C}}\mid  M^T=M , \Im M>0\}.
 $$
 %%%%%%%%%%%%%%%%%%%%%%%%%%%%%%%%%%%%%%%%%%%%%%%
 \begin{Prop}\label{uomegaest}
 Let $F\colon K\times\mathbf{C}^+ \to {\cal H}^n$ be a continuous map,
 $K$ compact, such that  
 for each point $u$ in $K$, its restriction to $\{u\}\times 
 {\mathbf{C}}^+$
 is holomorphic.
 Then for all $(u,\omega)$ in $K\times U_1$ where
 \begin{equation}\label{U1}
 U_1=\{\omega 
 \in 
 {\mathbf{C} }\mid   -1< x < 1, y \geq 
 1\},
 \end{equation}
 there are constants $C_{1}>0$ and $C_{2}$ such that
 \begin{equation}\label{ineq1}
 \det \Im F(u, \omega) \geq  C_{1} y^{-n} ,
 \end{equation}
and,  putting
 $F^{\prime}=\frac{\partial F}{\partial{\omega}}$,
\begin{equation}\label{ineq2} 
 \mid \det F^{\prime}(u,\omega)\mid  \leq  C_{2}.
\end{equation}
 \end{Prop}
 %%%%%%%%%%%%%%%%%%%%%%%%%%%%%%%%%%%%%%%%%%%%%%%%%%%%%%%%%
 \begin{pf}The inequality (\ref{ineq1}) would follow from
  Lemma 8.4 in \cite{RSZOKE95},
 but we include a proof here (along the lines of Sz\H{o}ke's) 
 for reference and because most of its ingredients  
 are needed to show inequality (\ref{ineq2}).  
 
 Inequality (\ref{ineq1}) follows from the compactness of $K$
  and the continuity of $F$ 
 simply by letting 
 $$C_{1}=\min_{{u\in K}} \det \Im F(u,\sqrt{-1}),
 $$ since
for fixed $u$ in $K$, along $\{ u\} \times U_{1}$
\begin{equation}\label{detimag}
 \det \Im F(u,\omega) \geq  (4 y )^{-n} \,
   \det \Im F(u,\sqrt{-1}).
\end{equation}     
To show (\ref{detimag}), recall from \cite{LESZ} that for fixed $u$ in $K$, the 
assumptions on $F$ imply that we have the Fatou representation
 \begin{equation}\label{Fharm}
\Im F(u,\omega)=y\,\alpha (u) + 
\frac{y}{\pi}\int_{-\infty}^\infty\frac{d\mu(u,t)}{(x -t)^{2}+y^{2}}
\end{equation}
for $\alpha(u)$ a non-negative symmetric constant matrix and  
$d\mu(u,t)$ 
a symmetric non-negative matrix of signed Borel measures on the real 
line, i.e.,  for all points $(v_1,\cdots ,v_n)$ in ${\mathbf{R}}^n$ 
\begin{equation}\label{measurenneg}
\sum_{ik}d\mu_{ik}(u,t) v_iv_k\geq 0
\end{equation}
as a measure,
where for all $1 \leq i,j\leq n$
$
\int \frac{\mid d\mu_{ij}(u,t)\mid }{1+t^2}<\infty .
$

Now , for $ |x| < 1$ and $y \geq 
 1$
we have that for all $t$ in ${\mathbf{R}}$,
$\frac{1}{1+t^2} < \frac{4 y^2}{(x-t)^2 + y^2}$
and thus, in virtue of (\ref{measurenneg}) and by the identity 
obtained  
from letting $y=1$ in (\ref{Fharm}), we get 
\begin{equation*}
\Im F(u,\omega) 
\geq y\,\alpha(u)  + 
\frac{1}{4 y\pi}\int_{-\infty}^\infty\frac{d\mu(u,t)}{t^2 +1}
=\frac{ (y^2-1)}{4y} +\frac{1}{4 y}\Im F(u,\sqrt{-1}) 
\end{equation*}
and so, 
\begin{equation}\label{FINEQ}
\Im F(u,\omega) 
\geq
\frac{1}{4y}\Im F(u,\sqrt{-1}) .
\end{equation} 
Finally, take determinants and recall that $M_1\geq M_2$ for $M_i$ 
non-negative real symmetric matrices implies $\det M_1 \geq \det M_2$. Thus, 
inequality (\ref{detimag}), and hence (\ref{ineq1}), are proved.

We now show inequality (\ref{ineq2}).  
By differentiating $F$ with respect to $\omega$, using (\ref{Fharm}) 
and that $F$ is holomorphic in $\omega$,  
\begin{equation}\label{Fprimeharm}
F^\prime (u,\omega)=\alpha (u)  + 
\frac{1}{\pi}\int_{-\infty}^\infty\frac{d\mu(u,t)}{(\omega -t )^2}
\end{equation}
and thus, for all $1\leq i,j\leq n$,   
 \begin{equation}\label{Fprimeharmij}
\mid F_{ij}^\prime (u, \omega) \mid 
\leq \mid \alpha_{ij}(u) \mid  + 
\frac{1}{\pi}
\int_{-\infty}^\infty
\frac{\mid d\mu_{ij}(u, t)\mid }{\big{(}(x-t )^2+y^2\big{)}^2}.
\end{equation}
%%%%%%%%%%
It is not hard to find appropriate bounds for the terms in the 
right-hand side above.
%%%%%%%%%%%%%%
Indeed, since $\alpha(u)$ is a non-negative symmetric matrix, then 
for all $1\leq i,j\leq n $
\begin{equation}\label{alphabd}
\alpha_{ii}(u)\geq0, \qquad \mid \alpha_{ij}(u)\mid ^2 \leq \alpha_{ii}(u) 
\alpha_{jj}(u).
\end{equation}
But, from (\ref{measurenneg}),  by putting $y=1$ in formula 
(\ref{Fharm})
we get for all $1\leq i \leq n$ 
$$ \alpha_{ii}(u) \leq  \Im F_{ii}(u, \sqrt{-1}),$$
and thus 
\begin{equation}\label{ivt0}
\mid \alpha_{ij}(u) \mid  \leq \max_{\{u \in K, \,\,1\leq k\leq n\}}
 \mid \Im F_{kk}(u, \sqrt{-1})\mid .
\end{equation}
%%%%%%%%%
Similarly, 
from (\ref{measurenneg}), for all $1\leq i,j \leq n$
\begin{equation}\label{ineqmeasures}
d\mu_{ii}(u,t)\geq 0, \quad \mid d\mu_{ij}(u,t)\mid ^2
\leq d\mu_{ii}(u,t) d\mu_{jj}(u,t),
\end{equation}
and thus from (\ref{ineqmeasures}) and the Cauchy-Schwarz inequality it follows, 
\begin{equation*}\label{Holder}
\Big{(}
\int_{-\infty}^{\infty}
\frac{\mid d\mu_{ij(u,t)}\mid}{\big{(}(x-t)^2 +y^2\big{)}^2}
\Big{)}^2
\leq 
\int_{-\infty}^{\infty}\frac{d\mu_{ii}(u,t)}{\big{(}(x-t)^2 +y^2\big{)}^2}
\int_{-\infty}^{\infty}\frac{d\mu_{jj}(u,t)}{\big{(}(x-t)^2 +y^2\big{)}^2}.
\end{equation*}
%%%%
Now,  for each of the factors in the right-hand side of the inequality above 
 we have by the inequalities on  the left in 
 (\ref{ineqmeasures}) and in (\ref{alphabd}),  
 always
for $(u, \omega )$ in $K\times U_1$,
\begin{equation*}
    \begin{align}\notag
\int_{-\infty}^{\infty}\frac{d\mu_{ii}(u,t)}{\big{(}(x-t)^2 +y^2\big{)}^2}
 & \leq
\int_{-\infty}^{\infty}\frac{d\mu_{ii}(u,t)}{(x-t)^2 +1 }
\leq \Im F_{ii}(u, x+\sqrt{-1})-\alpha_{ii}(u)
\\ \notag
&
\leq
\Im F_{ii}(u, x+\sqrt{-1}).\notag
\end{align}
\end{equation*}
%%%
So, 
\begin{equation}\label{ivt}
\Big{(}
\int_{-\infty}^{\infty}\frac{|d\mu_{ij(u,t)}|}{\big{(}(x-t)^2 +y^2\big{)}^2}
\Big{)}^2
\leq 
\Im F_{ii}(u, x+\sqrt{-1})
\Im F_{jj}(u, x+\sqrt{-1})
\end{equation}
and hence by (\ref{Fprimeharmij}),  (\ref{ivt0}) and (\ref{ivt}), 
\begin{equation}\label{estimateFprime}
|F^\prime_{ij}(u, \omega)|
\leq
c := 2 \max_{\{u \in K,\, ; \mid x \mid \leq 1 \,;1\leq k,l\leq n\}} 
\mid \Im F_{kl}(u, x+\sqrt{-1})\mid ,
\end{equation}
which implies
$$\mid \det F^\prime(u, \omega) \mid \leq  n! c^n , 
$$
showing the inequality (\ref{ineq2}). 
\end{pf}
%%%%%%%%%%%%%%%% End of Proof of Prop \ref{ESTIMATEM}}%%%%%%%%%%
\begin{pf*}{End of Proof of Proposition \ref{ESTIMATEM}}
To find the growth estimate for $\delta$ throughout $\{\tau \geq 
1\}$ (and hence on all of $X$) from 
the estimates along the leaves obtained earlier we 
simply use that the level set $\{ \tau =1 \}$ being compact   
can be covered by a finite number of images of geodesic flow boxes.
In fact, it is well-known that given $z$ in $X$ with $\tau(z)=1$ 
we may find positive numbers $r_{1}$ and $r_{2}$ and a real analytic map 
$$f \colon \{ u= (u_{1}, \cdots, 
u_{2n-2}) \in {\mathbf{R}}^{2n-2}| \sum u_{i}^{2}\leq r_{1} \} \times 
[-r_{2}, r_{2}]\to \{\tau=1 \}$$ with $f(0,0)=z$ 
and such that if $\gamma_{u}$ denotes the unit speed geodesic with initial 
conditions $\dot{\gamma}_{u}(0)= f(u,0)$ and 
$\gamma_{u}(0)=\pi(f(u,0))$ 
then, $\gamma_{u}(x)=f(u,x)$ for $(u,x)$ in the domain above.
(We already used these maps at the end of the proof of Prop. \ref{EXTENSION}
but here we are more explicit.)
Letting $u$ vary in ${\mathbf{ R}}^{2n-2}$ so that $\sum u_{i}^{2}\leq 
r_{1}$ we get a family
 of geodesics 
$\gamma_{u}$ and leaves ${\cal{L}}_{\gamma_{u}}$ that depend real 
analytically on $u$. Similarly for the vector fields, Jacobi and Fermi,  along   
$\gamma_{u}$, for their extensions to ${\cal{L}}_{\gamma_u}$ and for the 
parametrizations $P_{\gamma_u}$.
Performing the construction in section 4, now with $u$ as a real 
analytic parameter,  
we obtain a family of matrices $h_{\gamma_{u}}(\omega)$
of meromorphic functions as in (\ref{etahxidef}) depending real 
analytically on $u$. We remark that the poles of these matrices as 
functions in ${\mathbf{C}}$ lie on the 
``real axis'', hence away from the region where we use them for the 
estimates.  Thus 
the function 
$$
F(u, \omega):= h_{\gamma_u}(\omega)
$$
satisfies the conditions of Proposition \ref{uomegaest} (exceedingly 
so, since only continuity in $u$ would suffice there).
The images of finitely many of maps like the $f$ above will cover the 
level set $\{\tau=1\}$ and it follows that we can find a pair of 
constants $C_1$ and $C_2$
as in Proposition \ref{uomegaest} appropiate for all the 
corresponding functions $F's$. 
The estimate (\ref{ESTIMATEMf}), in fact a stronger one, 
\begin{equation}\label{strineq}
1/\delta \le c \tau^{n}, 
\end{equation}
 now follows from applying Proposition 
\ref{PDELTA}.
\end{pf*}
The Proof of Theorem 
1 is at this point complete. 
%%%%%%%%%%%%%%%%
%%%%%%%%%%%%%
%%%%%%%%%%%%%%%%%%%%%%%%%%%%%%%%%%%%%%%%%%%%%%%%%%%%%%%%%%%%%%%%%%%%
%%%%%%%%%%%%%%%%%%%%%%%%%%%%%%%%%%%%%%%%%%%%%%%%%%%%%%%%%%%%%%%%%%%
%%%%%%%%%%
\section{{\bf Rationality of the holomorphic metric}}\label{RATHMET}
%%%%%%%%%%%%%%%%%
%%%%%%%%%%%%%%%%%%%%%%%%%%%%%%%%%%%%%%%%%%%%%%%%%%%%%%%%%%%%%%%%%%%%%
%%%%%%%%%%%%%%%%%%%%%%%%%%%%%%%%%%%%%%%%%%%%%%%%%%%%%%%%%%%%%%%%%%%%
%%%%%%%%%%%%%%%
%%%%%%%%%%%%%%%%%
Since the Riemannian metric $g$ 
induced by $\sqrt{-1}\partial\Bar{\partial}\tau$ 
along ${M}$
is real analytic, 
it extends as a holomorphic object $\tilde{g}$ to a neighborhood 
${\cal U}$ of ${M}$ in ${X}$, i.e., 
$\tilde{g}$ is the holomorphic ``metric'' in $T^{1,0}{\cal U}$ 
continuing the 
${\mathbf{C}}$-bilinear extension to $T^{1,0}{X}\mid_{M}$ 
of the Riemannian metric on ${M}$ defined by
$$\tilde{g}(\xi^{1,0}, \eta^{1,0})=
g(\xi,\eta)$$
for all $\xi$ and $\eta$ in $T{M}$.
Along each leaf of the 
Monge-Amp\`{e}re foliation $\tilde{g}$ has the holomorphic Fermi 
fields 
as an ``orthonormal'' basis, i.e., $\tilde{g}=\sum_{i}^{n}
F_{i}\otimes F_{i}$.
It follows from an argument similar to the one used in our proof of 
the extension of the volume form in Proposition \ref{EXTENSION}
that $\tilde{g}$ extends to a holomorphic non-degenerate 
section of the symmetric 
product of the holomorphic cotangent bundle 
on the entire tube ${X}$.

We will show that such a section is actually rational on the 
affine ${X}\subset C^{N}$. Equivalently, if $\bar{X}$ is a smooth
compactification of $X$, i.e., $X \subset \bar{X}$ and $D := \bar{X}
\setminus X$ is a subvariety of dimension less than $n = {\mbox{dim}}
\, X$, then $\tilde{g}$ has a meromorphic extension to all of $\bar{X}$
which is holomorphic and non-degenerate on $X$.

\begin{thm}\label{grational}
Let $\big{(}{X}, \tau\big{)}$ be an entire Grauert tube.
The $\tau$-induced Riemannian metric on $\{\tau=0\}={M}$ 
extends as a rational 
holomorphic metric on all of $\tilde{{X}}$.
Consequently, the holomorphic volume form in Proposition 
\ref{EXTENSION}
extends as a rational
form on $\tilde{X}$.
\end{thm}

This will follow from Proposition \ref{length}
below, whose proof involves
our previous estimates and some facts from 
\cite{DEMAILLY}. Especially  Theorem 8.5 there,  which shows 
that the transcendence degree over ${\mathbf C}$ 
of the field of fractions of $\phi$-polynomials, $K_{\phi}$, 
is at most 
$n =$ dim ${X}$, implies 
that any $\phi$-polynomial is a rational 
function in ${X}$.

Let's put $\tilde{g}^{*}$
for the  dual metric on $(T^{1,0}{X})^{*}$. 
Recall that $\beta=\sqrt{-1}\partial\Bar{\partial}(e^{\phi})$ with 
$\phi$ as in Theorem 1.

\begin{Prop}\label{length}
The length of $\tilde{g}^{*}$ in the metric induced by $\beta$ satisfies
$$\parallel \tilde{g}^{*}\parallel _{\beta} 
\leq c_{1} e^{c_{2} \sqrt{\tau}}
$$
for constants $c_{i}$.
\end{Prop}
\begin{pf}
Along a parametrized leaf 
we have, outside the real discrete set $S$, omitting variables,  
$ F_{i} = 
\sum_{k=1}^{n} C_{ki} \xi_{k}^{1,0}
$, where $C=A^{-1}$.
So, 
$$\partial \Bar{\partial}\tau
  ( F_{i}, \Bar{F}_{j} ) = 
\sum_{k,l=1}^{n} C_{ki}\Bar{C}_{lj}
\partial \Bar{\partial}\tau \big{(}
 \xi_{k}^{1,0}, \xi^{0,1}_{l}\big{)}
$$ and thus, since by (\ref{ATAHPRIME})
$$
\quad \sum_{k=1}^{n}C_{ik}C_{jk}=h^{\prime}_{ij}, 
$$ we have, by (\ref{MIJD}) and (\ref{MIJ}), 
\begin{equation*}
\begin{align}\notag
 \sum_{i,j=1}^{n}
 \big{(}
 \partial \Bar{\partial}\tau ( F_{i}, \Bar{F}_{j}) \big{)}^{2} = &
 \; 
 y^{2}\sum_{i,j, k, l , s, t  =1}^{n} C_{ki}\Bar{C}_{lj} C_{si}\Bar{C}_{tj} 
 (\Im h)^{-1}_{kl} (\Im h)^{-1}_{st} \\ \notag
 =& \;
 y^{2}\sum_{i,j, k,l, s, t =1}^{n} C_{ki}C_{si} \Bar{C}_{lj} \Bar{C}_{tj}
 (\Im h)^{-1}_{kl} (\Im h)^{-1}_{st}\\ \notag
 =&\; y^{2}
  \sum_{k,l, s, t =1}^{n} h^{\prime}_{ks} \Bar{h}^{\prime}_{lt}
 (\Im h)^{-1}_{kl} (\Im h)^{-1}_{st}.\notag 
 \end{align}
\end{equation*}
It follows, since $\Im h$ is positive-definite and hence so is $(\Im h)^{-1}$,
that
$$
\sum_{i,j=1}^{n}
 \big{(}
 \partial \Bar{\partial}\tau ( F_{i}, \Bar{F}_{j}) \big{)}^{2} \leq 
y^{2}
  \sum_{k,l, s, t =1}^{n} \mid h^{\prime}_{ks}\mid 
  \mid \Bar{h}^{\prime}_{lt}\mid
 \big{(}(\Im h)_{kk} (\Im h)_{ll} (\Im h)_{ss}(\Im 
 h)_{tt}\big{)}^{-\frac{1}{2}}.
$$
Now we use Proposition \ref{uomegaest} with $F=h(u, \omega)$ to 
estimate the right-hand side above. 
From (\ref{FINEQ})
it follows that $$(\Im F)^{-1}(u, \omega)\leq 4 y (\Im F)^{-1}(u, 
\sqrt{-1})$$ (for $(u, \omega) $ in an appropriate $D_{r_{1}}\times D_{r_{2}}$)
and this together with the estimate (\ref{estimateFprime})
implies, after arguing as in the end of the proof of Proposition 
\ref{ESTIMATEM}, 
that for $\tau \geq 1$ there is a constant $c>0$ so that
\begin{equation}\label{taufermi}
\mid
 \partial \Bar{\partial}\tau \big{(} F_{i}, \Bar{F}_{j} \big{)} 
\mid
 \leq c \tau,
\end{equation}
where the inequality (\ref{taufermi}) with the same constant $c$ is valid 
for any set of holomorphic Fermi fields along any parametrized geodesic.  

%%%%%%%%%%%%%%%%%%%%%%%%%%%%%%%%%%%

Using (\ref{taufermi}) and  the expression 
for the $\beta$-metric 
(\ref{betametric}) it is now easy to verify
that there is a constant $c_{1}>0$ so that
\begin{equation}\label{betafermi}
\| F_i \|_{\beta}^2
\leq c_1 e^{2\, \sqrt{\tau}}.
\end{equation}
(For this calculation one recalls the splitting  induced by the kernel of 
$\partial \tau$ and the evaluation of the $F_{k}$ on  $\partial \tau$ 
 and $\partial \Bar{\partial} \tau$ as discussed in the proof of Prop. 
 \ref{EXTENSION}.)
Now, from  11.5(c) in \cite{DEMAILLY}, 
and with $C_1 \gg 0$, $\rho \equiv 0$, for any point $z
\in {X}$, we can find $f_1,\cdots ,
f_n$, such that $$f_i \in {\mbox{A}}_{\phi}^b$$
 (and
therefore ``$\phi$-polynomials'' by Lemma 11.3 in \cite{DEMAILLY}),  where 
$${\mbox{A}}^{b}_{\phi}
= \{ f\colon X\to {\mathbf{C}} 
{\textit{ holomorphic }}\;
\mid \; \exists C \geq 0 \;\mid \; \int_{X} \mid f \mid ^{b} e^{-C\phi} 
\beta^{n}< \infty
n\},
$$satisfying
$$df_1 \wedge df_2
\cdots  \wedge df_n \neq 0,$$ at $z \in {X}$. Here $b = \frac{2c}{1+c} \in
(0,2)$. Thus, the differentials $df_i$ give a rational frame for the
holomorphic cotangent bundle $T_{(1,0)}^{*}(X)$ on a Zariski open set
of $X$. The coefficients of the dual holomorphic metric $\tilde{g}^{*}$ on
$T_{(1,0)}^{*}(X)$ are given, in this frame and on this Zariski open
set, by
$$\tilde{g}^{*}( df_i, df_j),$$ and we have succeeded in showing that $g$ is
rational if these coefficients are ``$\phi$-polynomials'' as indicated 
earlier. 
In fact, 
we can
 show that they are in ${\mbox{A}}_{\phi}^{\frac{b}{2}}$. To see
this latter, consider
$$\int_{X} \| \tilde{g}^{*}(df_i, df_j) \|_{\beta}^{\frac{b}{2}} e^{-\tilde{C}
\phi} \beta^n.$$ The integrand can be estimated: 
$$ \| \tilde{g}^{*}(df_i,
df_j) \|_{\beta}^{\frac{b}{2}} \leq \|\tilde{g}^{*}\|_{\beta}^{\frac{b}{2}} \; \|
df_i \|_{\beta}^{\frac{b}{2}} \; \| df_j \|_{\beta}^{\frac{b}{2}}.$$
But by (\ref{betafermi}), we can estimate $$\| \tilde{g}^{*}
\|^{\frac{b}{2}} \leq c_1' e^{c_2' \; \sqrt{\tau}}$$ which can be absorbed into 
the weight term
if $\tilde{C} \gg 0$. Again, perhaps taking $\tilde{C}$ larger, we can
then use the Cauchy-Schwarz inequality to conclude that the
coefficients are in ${\mbox{A}}_{\phi}^{\frac{b}{2}}$. This last step
uses 11.3 (b) in \cite{DEMAILLY}.
This shows Proposition \ref{length} and, as indicated above, in virtue 
of Theorem 8.5 in \cite{DEMAILLY}, Theorem \ref{grational} is thus proved.
\end{pf}

The rationality of $\tilde{g}$ is of interest 
in connection with the problem of characterizing 
those Riemannian manifolds that appear as ``centers'' of entire 
Grauert tubes,  a question we hope to come back to 
in the future. But for now, we will 
derive from Theorem 3 some  
natural corollaries related to that problem.

In general, an entire Grauert tube 
$({X}, \tau_1)$ is non-rigid in the sense
that it might support another Monge-Amp\`{e}re exhaustion 
$\tau_{2}$ different from a constant multiple
of $\tau_{1}$ with the same zero set ${M}$. 
Thus the Riemannian structures induced in ${M}$ by the different
Monge-Amp\`{e}re  functions 
are non-isometric, even non-homothetic.
(Note that such $\tau_{2}$ is necessarily
unbounded. This is so because the maximum on the level set $\{\tau_{1}= r\}$ 
of any plurisubharmonic function such as  
$\tau_{2}$ is increasing convex as a function of $r$. See Corollary 
6.6 in \cite{DEMAILLY}.)

The rationality theorem above yields, under an additional assumption
on $D = \bar{X} \setminus {X}$, where $\bar{X}$ is a smooth,
projective compactification of the affine manifold $X$, the
``rigidity'' of the Riemannian volume in ${M}$ as follows.

\begin{Cor}\label{EQVOL}
Let $\tau_{i}$ be two Monge-Amp\`{e}re exhaustions on ${X}$
so that $({X}, \tau_{i})$ is an entire Grauert tube with
${M}=\{ \tau_{1}=0\}=\{\tau_{2}=0\}$.
If the divisor $D=\tilde{X}\setminus {X}$ is irreducible then
the Riemannian volume forms on ${M}$ corresponding to  
the Riemannian metrics $g_{i}$  induced by $\tau_{i}$ are, 
up to homothety,  identical.
\end{Cor}
\begin{pf}
The argument in the previous proof shows that the 
holomorphic forms ${\cal V}_{i}$ extending the Riemannian 
volume forms of $({M}, g_{i})$ 
extend in turn as meromorphic 
sections ${\frak V}_{i}$
of the canonical bundle $T^{n,0}(\tilde{{X}})$ of
the compactification $\tilde{{X}}$. By Proposition 
(\ref{EXTENSION})  ${\frak V}_{i}$ are non-vanishing on ${X}$. 
Thus 
$${\frak V}_{1}=f\,{\frak V}_{2}$$
where 
$f$ is a 
meromorphic function on $\tilde{{X}}$, 
holomorphic and non-vanishing on 
${X}$. Assume $f$ is not constant. Then if $D$ is the polar set of 
$f$ (or of $1/f$) by 
the Riemann Extension Theorem  $1/f$ (or $f$) extends as a 
non-constant holomorphic function  on 
the compact $\tilde{{X}}$, which is not possible.
\end{pf}

Examples of unbounded Grauert tubes with 
$D=\tilde{X}\setminus {X}$ an irreducible hypersurface are the
standard complexifications 
$G_{\mathbf C}/H_{\mathbf C}$
of compact symmetric spaces of rank one, 
${M}=G/H$, $G$ a compact connected semi-simple Lie group
and $H\subset G$ a closed subgroup. The Stein manifold 
$G_{\mathbf C}/H_{\mathbf C}$ 
can be real-analytically 
identified with the cotangent bundle of 
$G/H$, 
$T^{*}(G/H)$, in a $G$-equivariant way and so that 
$G/H\subset G_{\mathbf C}/H_{\mathbf C}$ is identified with  
the zero section. 
The Monge-Amp\`{e}re exhaustion 
$\tau$ is the function on $G_{\mathbf C}/H_{\mathbf C}$ 
corresponding to the length squared function on $T^{*}(G/H)$
induced by the metric on $G/H$. 
From the explicit realization of
$G_{\mathbf C}/H_{\mathbf C}$
as affine manifolds (cf. \cite{GPPMW})
we see that the hypothesis on the divisor at infinity 
for their compactification is satisfied.

Finally, note that the anti-holomorphic involution $\sigma$ operates
on holomorphic functions on $X$ by $\sigma(f)(z) =
\overline{f(\sigma(z))}$. This action leaves invariant the various
$\phi$-polynomial norms. Given the affine embedding of $X$ by
$\phi$-polynomials, $(f_1,...,f_N)$, we may expand to the affine
embedding with components $(f_1,...,f_n,\sigma(f_1),...,\sigma(f_N))$,
and after composing with a linear change of variable in ${\mathbf
C}^{2N}$, we may arrange that we have an affine embedding of $X$ into
${\mathbf C}^{2N}$ which is $\sigma$-equivariant, that is:

$$(F_1(\sigma(z)),...,F_{2N}(\sigma(z))) =
(\overline{F_1(z)},...,\overline{F_{2N}(z)}).$$  This will mean that the
projective closure $\overline{F(X)}$ of $F(X) \subset {\mathbf C}^{2N}$ inside
${\mathbf{P}}^{2N}$ will be a real variety invariant under complex
conjugation. By \cite{BM}, there is a smooth resolution $\widehat{F(X)}$
of $\tilde{X}$ with an
exceptional divisor with normal crossings which is invariant under the
lift of $\sigma$ to $\hat{\sigma}: \widehat{F(X)} \rightarrow \widehat{F(X)}.$
It is also easy to check that the rational metric $\tilde{g}$, for example,
is $\sigma$-invariant on $X$ in the natural sense.

%%%%%%%%%%%%%%%%%%%%%%%%%%%%%%%%%%%%%%%%%%%%%%%%%%%%%%%%%%
 %%%%%%%%%%%%%%%%%%%%%%%%%%%%%%%%%%%%%%%%%%%%%%%%%%%%%%%%%%
%%%%%%%%%%%%%%%%%%%%%%%%%%%%%%%%%%

\section{{\bf Proof of Theorem 2}}

%%%%%%%%%%%%%%%%%%%%%%%%%%%%%%%%%%%%%%%%%
%%%%%%%%%%%%%%%%%%%%%%%%%%%%%%%%%%%%%%%%%%%%%%%%%%%%%%%%%%%%%%
We continue identifying  the entire Grauert tube $({X}, \tau)$ 
with the tangent bundle $T{M}$ of ${M}=\{\tau =0\}$, 
endowed with the adapted complex structure ${\mathbf{J}}$. 
According to Lempert and Sz\H{o}ke 
(Theorem 4.1 in \cite{LESZ}) 
the metric $g$ on ${M}$ has non-negative Gauss curvature. By Gauss-Bonnet, 
the Euler 
characteristic of ${M}$ is non-negative, and so
the center of the tube is {\em isometric} to a flat two-torus $T^{2}$ or 
or flat Klein bottle \cite{LESZ},  
or diffeomorphic to a two-sphere $S^{2}$ or real projective plane.

By the uniqueness of the adapted
complex structure, in the case of the flat two-torus  the Grauert tube 
is biholomorphic to 
$\mathbf{C}^{2}\big{/} {\Lambda}\simeq \mathbf{C}^{*}\times
\mathbf{C}^{*}$, where $\Lambda$ is a lattice in $\mathbf{R}^{2}$, 
 or a ${\mathbf Z_2}$ quotient of it, in the case of the Klein bottle.

Again by functoriality and uniqueness of the adapted complex structure
the projective plane case is subsumed into the one
which we now treat, ${M}$ diffeomorphic to $S^{2}$.

By Theorem 1 the Grauert tube ${X}$
can be imbedded algebraically in some ${\mathbf C}^{N}$ so that 
${X}$ is a Zariski open set of a smooth projective 
compact variety $\widetilde{{X}}$. Thus, $\tilde{{X}}$ 
is a  smooth projective compactification 
of ${X}$.
The ``divisor at $\infty$'' is given by, 
$$\widetilde{{X}}\setminus{X}=D= \bigcup_{i \in I} D_{i}
$$
where $D_{i}$ are Riemann surfaces which, by blowing-up if necessary, 
 will be assumed to intersect 
normally, i.e.,
$
\textit{ if } i\neq j, \quad D_{i}\cap D_{j} =\{ \textit{ point }\} \textit{ or } =\emptyset
$
and 
$
 \textit{ if }
i\neq j \neq k \quad 
D_{i}\cap D_{j}\cap D_{k}= \emptyset.$
Here each $D_i$ is an irreducible component obtained as the closure
of a connected component of the regular part of $D$.

We now argue as Ramanujam in \cite{RAMANUJAM}.
In what follows homology and cohomolgy is taken with rational 
coefficients $\mathbf{Q}$, unless indicated differently.
Let $N$ be a tubular neighborhood of 
$D\subset \tilde{{X}}$, a four-dimensional manifold 
with boundary $\partial N$.
For  $r>0$,
 $$\partial{N} \simeq \{ \tau = r\} \simeq {\mathbf{SO}}(3, 
 {\mathbf{R}})\simeq S^{3}/{\mathbf{Z}}_{2}.$$ 
Thus $\partial{N}$ is connected and  so is $D$, a deformation retract 
of $N$. Moreover
$
H^{1}(\partial{N})\simeq H_{2}(\partial{N})=0$.
On the other hand
$$H^{k}(N, \partial N)\simeq H_{4-k}(N) \simeq H_{4-k}(D)$$ for $0\leq k 
\leq 4,
$ the first identification by the relative Poincar\`{e} duality.
But, since 
$H_{3}(D)=0$, $H^{1}(N, \partial N)=0$, and hence  
it follows from the exact 
sequence 
$
\ldots  \mapsto H^{1}(N, \partial{N}) \mapsto
H^{1}(N)\mapsto H^{1}(\partial{N}) \mapsto \ldots
$, that $H^{1}(N)=0$ and so 
$
H^{1}(D)=0.
$
As a consequence, each $D_{i}\simeq \mathbf{P}^{1}$, 
the complex projective plane, and moreover there 
are no ``loops'', i.e.,  there is no sequence 
of distinct values  $\{ i_{1}, \dots , i_{N} \} $
for which
$D_{i_{1}}\cap D_{i_{N}}\neq \emptyset$ and 
$
D_{i_{k}}\cap D_{i_{k+1}}\neq \emptyset, \textit{ for all }
 1\leq k\leq N-1.
$

The {\em intersection matrix} for $D$ has  
coefficients $(I_{ij})$ 
defined by:
$I_{ii}=$ self-intersection number of $D_i$, 
$I_{ij} =1$ if $D_{i}\cap D_{j}\neq 
 \emptyset $ and $i\neq j$ and 
 $I_{ij}=0 $ otherwise. From the discussion above, 
 to the data $\{ D_{i}\}$ there is associated  
a weighted tree $T$ with a 
vertex $v_{i}$ for each  curve $D_{i}$ with
the weight on $v_{i}$ equal to $I_{ii}$, two
distinct vertices $v_{i}$ and $v_{j}$ being linked if and only if
$D_{i}$ and $D_{j}$ intersect.
 It is well-known that the fundamental group of  
$\pi(\partial N)$ 
is isomorphic
to the ``fundamental group'' $\pi(T)$ of  the tree $T$,  
a  presentation of which is given by the set of vertices
$\{v_{i}\}$ subject to the 
set of relations $\{ R_{i}, i\in I\}$
\begin{equation}\label{present0}
\{ R_{i} \}=
\begin{cases}
v_{i} v_{j}^{I_{ij}} = v_{j}^{I_{ij}} v_{i} \quad 1\leq i,j \leq N 
\\
1=v_{1}^{I_{i1}}  \cdots  v_{N}^{I_{iN}} \quad 1\leq i\leq N 
\end{cases}
\end{equation}
where in the second line in the relations $\{ R_{i}\}$ the product of the 
$v_{j}^{I_{ij}}$ a fixed order (irrelevant up to isomorphism)
of all the vertices is assumed.

A classical theorem of Mumford 
\cite{MUM} asserts that 
if the fundamental group of the connected tree associated to a
connected divisor with negative-definite intersection matrix is
trivial, then the divisor can be blown-down to a smooth point.  
We will later refer to this result as ``Mumford Theorem''.
In the next propositions
we determine how to choose subtrees of $T$ to which we can 
apply Mumford's theorem.

As explained after the statement of Corollary \ref{EQVOL}, we can 
assume the resolution $\tilde{X}$
to be so that the canonical anti-holomorphic involution
$\sigma \colon {X} \to {X}$ extends as an 
anti-holomorphic involution 
$
\tilde{\sigma}\colon 
\tilde{{X}}
\to 
\tilde{{X}}$, 
which then acts on  $D$.
Since $\tilde{{X}}$ is a complex compact two-dimensional
  K\"{a}hler manifold, 
the cup product restricted to 
$H^{1,1}(\tilde
{{X}}) \cap H^{2}(\tilde{{X}}, {\mathbf Z})$ 
has exactly
a 1 dimensional subspace where it is 
positive definite, by the Hodge
index theorem \cite{GH}.
We will next show that such a subspace is, after bowing down,
 the span of one of the $D_{i}$'s. 
So the strategy now is first to identify the part 
of $D=\tilde{{X}}\setminus{X}$ that is invariant by 
$\tilde{\sigma}$, which is done in Proposition \ref{FIX}
and then to show that its transform after blowing down the complement, 
has positive self-intersection, which is proved in Proposition 
\ref{BD}. The final step will be to show that the surface is actually 
rational.

%%%%%%%%%%%%%
\begin{Prop}\label{FIX}
We may assume $D=D_{0} \cup_{j>0}\{ D_{j}\cup D_{-j} \}$
with $D_{j}\neq D_{k}$ for $j\neq k$, 
$\tilde{\sigma} D_{0}=D_{0}$ and $ \tilde{\sigma} D_{\pm j}=D_{\mp j}
$, 
with $\tilde{\sigma}|_{D_{0}} \neq {\textit{ Identity  }}$.
\end{Prop}
%%%%%%%%%%%%%
\begin{pf}
Since the compactification  is assumed  invariant
with respect to $\tilde{\sigma}$, we have that
$\tilde{\sigma}D=D$. To find out how many components  $D_{i}$  
 are non-trivially permuted by $\tilde{\sigma}$, i.e., 
 $\tilde{\sigma}D_{i}\neq D_{i}$, 
 we consider the action of $\tilde{\sigma}$ on
 $H^{*}(\tilde{X})$.

Because $D$ has the homotopy type of a wedge product
of $N$ $2$-spheres, 
$H^{2}(D)\simeq {\mathbf Q}^N.
$
%%%
But,
$$H^{k}(\tilde{{X}}, D) 
\simeq H^{k}_{cpt.}(X)\simeq H_{4-k}(X)
\simeq H_{4-k}( S^{2}),$$
where $H_{cpt.}$ denotes cohomology with compact support, 
and by  
the exact sequence $
\ldots \mapsto H^{k}(\tilde{{X}}, D) \mapsto
H^{k}(\tilde{{X}})\mapsto H^{k}(D)\mapsto  \ldots
$
we get 
\begin{equation}\label{h13}
H^2(\tilde{{X}})\simeq {\mathbf Q}^{N+1}, \quad 
H^{k}(\tilde{\cal{X}})\simeq 0, \quad k=1 \textit{ or } k=3.
\end{equation} 
Since $\tilde{{X}}$ is connected, orientable,
and the complex dimension of $\tilde{{X}}$ is even, 
$\tilde{\sigma}_{*}$ is the identity on $H^{0}(\tilde{{\cal 
X}})\simeq {\mathbf Q}$
and also on $H^{4}(\tilde{{X}})\simeq {\mathbf Q}$. 
Now, 
$\tilde{\sigma}$ restricted to $\{ \tau = 0\}$ is the identity map, 
and this gives an additional 
1-dimensional subspace of 
$H^2(\tilde{{X}})$ where $\tilde{\sigma}_*=I$.
But each connected component of the fixed-point set 
of the anti-holomorphic involution 
$\tilde{\sigma}$ is a totally real 
submanifold of $\tilde{{X}}$. It follows that 
if $\tilde{\sigma} D_i=D_i$
then $\tilde{\sigma}$ can  not be the identity on such 
$D_{i}$ but  
must induce an antiholomorphic involution. In particular
$\tilde{\sigma}$ restricted to any invariant $D_i$  is orientation reversing 
and induces $-$Identity on the 
subspace 
in $H^2(\tilde{{X}})$  spanned by such  $D_i$.

So, from the discussion above, letting
$\nu$ be the number of the $D_{i}$'s 
invariant by $\tilde{\gamma}$, 
by the Lefschetz Fixed Point Theorem  
$$\sum_{k=0}^4 (-1)^k \textit {trace } ( \tilde{\sigma}_* 
|_{H^k(\tilde{X})} )
= 1-0+(1-\nu)-0 +1 =\epsilon (S^2) = 2,$$
for the total 
contribution to the trace 
of those $N-\nu$ subspaces of $H^2 (\tilde{{X}})$ corresponding to the
$D_i$  that are non-trivially permuted
is zero. 

Thus, 
$\nu =1$,  i.e., there is only one curve of $D$ left invariant    
by $\tilde{\sigma}$ which we will denote by $D_0$ 
and call the ``central curve''.
\end{pf}
%%%%%%%%%%%%%%
\begin{Prop}\label{euler}
Let $\tilde{{X}}$ and  $D$ as above, with $N$ the number of 
irreducible components of $D$. Then the Euler and holomorphic Euler 
characteristics of $\tilde{X}$ are
$$\epsilon(\tilde{{X}})=3+N, \quad 
\chi (\tilde{{X}}, {\cal O}_{\tilde{{X}}})=\sum_{p=o}^2 (-1)^p h^{0, p}=1.$$
\end{Prop}   
%%%%%%%%%%%%%%%
\begin{pf}
For the Euler characteristic  $\epsilon(\tilde{X})$ use (\ref{h13}) 
and for the holomorphic 
Euler characteristic, the compactness of $\tilde{{X}}$ and 
Hodge decomposition. Here the simple connectedness 
of
$\tilde{{X}}$ implies  $h^{01}=\frac{1}{2}h^{1}=0$. Also the  subspace of  
$H^{1, 1}(\tilde{{X}})$
generated by $D$ has 
codimension 1 
in $H^2(\tilde{{X}})$ and thus $h^{0,2}=0$.
\end{pf}
%%%%%%%%%
\begin{Prop}\label{BD}
The set of curves $\bigcup_{i\neq 0}D_{i}$ 
can be blown-down.
\end{Prop}
%%%%%%%%%%%%%%
\begin{pf}
Let $F_{\pm k}$ be the connected components of
$$\bigcup_{i\neq 0} D_{i}=\bigcup_{k>0}^{N^{\prime}} F_{k}\cup 
F_{-k},$$ where 
$2N^{\prime}$ is the number of 
$D_{i}$'s meeting the central curve. Here 
$\tilde{\sigma}F_{k}=F_{-k}$,  and so, 
$F_{k}$ being orthogonal to $F_{-k}$ and the intersection product 
invariant by $\tilde{\sigma}$, by 
the Hodge-Riemann index theorem,  
the intersection product on  each $F_{k}$ must be negative definite.
Thus, by Mumford Theorem,
a connected component $F_{k}$  
can be blown-down provided that   
the fundamental group of the branch $T_{k}$ on $T$ 
associated to $F_{k}$ satisfies 
\begin{equation}\label{ptt}
\pi(T_{k})=1.
\end{equation}
Now, consider the group $\pi(T)\big{/}(v_{0})$ which  may be presented
by adding the relation $v_{0}=1$ in (\ref{present0}) or, 
 equivalently,  by the  set of generators 
 $\{ v_{1}, \cdots , v_{N}\}$ with the relations 
$$
\tilde{R}_{i}=
\begin{cases}
 \{ R_{i} \}
  \textit{ for } i\neq 0 
\\ 
1=v_{j_{1}}\cdots v_{j_{2N^{\prime}}}
\end{cases}
$$
where in the last relation  
the generators involved are those 
corresponding to the vertices linked to $v_{0}$ in the original tree 
$T$. Except this last relation, there are no
further relations among generators corresponding to vertices in 
different  branches 
$T_{k}$'s. Thus,  
this group may be identified with
the amalgamated free product
\begin{equation}\label{afp} (\pi(T_{-N^{\prime}}) * \cdots * \pi(T_{-1}) )*\pi(T_{1}) )
*\cdots *\pi(T_{N^{\prime}}) )
\big{/}
(v_{j_{1}} * \cdots * v_{j_{ 2 N^{\prime} } })
.\end{equation}
We also need the following 
\begin{Prop}\label{amalgam}\cite{MUM}
Let $G_{i}$, $i=1,2,3$ be any non-trivial groups, and let 
$a_{i}\in G_{i}$ be 
arbitrary elements.
Then the amalgamated free product 
$$\big{(} G_{1}*G_{2}*G_{3}\big{)}
\big{/}\big{(}a_{1}*a_{2}*a_{3}\big{)}$$
is non-trivial non-cyclic.
\end{Prop}
Since $\pi(T)\big{/}(v_{0})$ is either trivial or 
${\mathbf{Z}}_{2}$ 
from (\ref{afp}),  Proposition \ref{amalgam} and Mumford Theorem
applied inductively, we may assume that, after blowning down pairs of 
branches $\{ F_{-k}, F_{k} \}$ in $D$,  there are at most 
two branches $F_{-}$ and $F_{+}$ say, interchanged by the 
involution $\tilde{\sigma}$ with  $\pi(T_{+}) \simeq \pi(T_{-})$. In 
particular, 
we have
$$\pi(T)\big{/}(v_{0}) \simeq \big{(} \pi (T_{+}) * \pi (T_{-}) 
\big{)}
\big{/}\big{(} v_{+}*v_{-} \big{)}.$$
\vskip3mm
Claim:{\em \;\;$\pi(T_{+}) \simeq \pi(T_{-})\simeq 1$.}

{\em Proof of Claim.} Assume $T_{+}$ has $K$ vertices. 
In an appropiate basis we 
write down the intersection matrix restricted to the span of the 
divisor $D$ as the matrix with integer 
coefficients, 
%%%%%
$${\frak{I}}=
\begin{pmatrix}
 a    & C & C \\
 C^{T}& B & 0\\
 C^{T} & 0 & B
\end{pmatrix} \in {\mathbf{Z}}^{(K+1)^{2}}
$$   
%%%% 
where, the superscript $''T''$ as earlier indicates transpose,  
each block 
$B$ is a symmetric matrix in ${\mathbf{Z}}^{K\times K}$ corresponding 
to the $K$ vertices in $T_{+}$ or to the $K$ vertices in $T_{-}$, 
$C$ is a $1\times K$ matrix with one entry equal to $1$ and all the 
other ones equal to $0$,  and finally 
$a=I_{00}$ is the weight of the central vertex
$v_{0}$. Since the intersection form is positive-definite in 
the span of the 
central curve and negative-definite in the branches $F_{+}$ and 
$F_{-}$, we have
$a>0$, and rank $B=K$.

To prove the Claim  all we need to show 
is that $a>1$. Indeed, if $a>1$, then take
the groups $G_{+}=\pi(T_{+})$, $G_{-}=\pi(T_{-})$ and 
$<v_{0}>$, the
cyclic group of order $a$ generated by $v_{0}$. Let $v_{\pm}$ the vertex in 
$T_{\pm}$ linked to $v_{0}$. Then, arguing as in  the discussion preceding 
expression (\ref{afp}) we see that 
\begin{equation}\label{pitree}
{\mathbf{Z}}_{2}¥
\simeq \pi(T)=\big{(}\pi(T_{+}) * \pi(T_{-}) * <v_{0}> \big{)} \big{/} 
(v_{+}*v_{-}*v_{0}\big{)}
\end{equation}
which,  unless $\pi(T_{\pm})\simeq 1$, 
would contradict Proposition \ref{amalgam}.

So, assuming $a=1$,
 we will get to a contradiction.
Consider the matrices
$$
U_{1} = \begin{pmatrix}
1&-C& -C\\
C^{T}&I&0\\
C^{T}&0&I
\end{pmatrix}, \quad 
U_{2} = \begin{pmatrix}
         1& 0 & 0 \\
         0 & I & 0 \\
         0 & I & I
\end{pmatrix}, \quad 
U_{3} = \begin{pmatrix}
         1& 0 & 0 \\
         0 & I &-I \\
         0 & 0 & I
\end{pmatrix}
$$
$
\in {\mathbf{Z}}^{(K+1)^{2}}
$ with $\det U_{i}=1$ for $i=1,2,3$.
Then, 
$$ U_{3}U_{1}^{T}{\frak{I}}U_{1}U_{2} = 
 \begin{pmatrix}
         1 & 0              & 0 \\
         0 & B - 2 \;C^{T}C & -C^{T}C\\
         0 & 0              &    B
\end{pmatrix}
$$
Thus, since $\pi(T)={\mathbf{Z}}_{2}$, 
 $$2= |\det{\frak{I}}|= |\det B | |\det (B-2\;C^{T}C)|.$$
 On the other hand,  since we are assuming $a=1$,  
 it follows from (\ref{pitree}) 
that $|\det B |>1$.  
Also,  $B-2 C^{T}C$ has integer coefficients, so
we must have 
\begin{equation}\label{deteq}
|\det B|=2,\quad |\det (B-2\;C^{T}C)|=1.\end{equation}
But if $|\det B|=2$ then 
$2 B^{-1}$ is well defined in ${\mathbf{Z}}^{(K+1)^{2}}$, 
$I-2B^{-1}C^TC$ has integer coefficients and hence,
\begin{equation}\label{detint}
\det (I-2B^{-1}C^TC)\in \mathbf{Z}.\end{equation}
Then, by the right handside equality in (\ref{deteq}) 
$$1=|\det (B-2\;C^{T}C)|=|\det B  | |\det (I-2\;B^{-1}C^{T}C)|,$$
a contradiction in virtue of (\ref{detint}) and the left handside equality 
in (\ref{deteq}).

We showed $a>1$ and $\# \pi(T_{\pm})=|\det{B}|=1$, thus proving 
the Claim. This finishes the proof of 
Proposition \ref{BD}.
\end{pf}

\subsubsection{Rationality of the surface}
An algebraic surface is called rational if it is 
birational to the complex projective plane ${\mathbf P}^2$, which 
is equivalent to saying that it 
may be obtained from it 
by blowing up points  followed by blowing down curves. 

Because of the results in the previous sections , the divisor at infinity in the 
algebraic surface $\tilde{X}$, $\tilde{X}\setminus X$, may be 
assumed to consist of just a $\mathbf{ P}^{1}$ with self-intersection 
$a=2$. This together with the topology of $X$, we will show, 
implies that $\tilde{X}$ verifies the conditions of the 
classical Castelnuovo Criterion for Rationality (cf. \cite{GH}), 
namely zero first Betti number 
and zero geometric genus $P_2$ (\ref{gg}). 
Moreover, the topology of $\tilde{X}$ will  force it to be ``minimal'', 
i.e. to
contain no curve with 
self-intersection $-1$, and thus to belong to a certain family of 
rational surfaces ${\mathbf F}_k$ which are projective bundles over 
${\mathbf P}^1$.
From here it will follow that $\tilde{{X}}$ is 
$F_{0}={\mathbf P}^1 \times {\mathbf P}^1$.

\begin{Prop}\label{MINRATIO}
$\tilde{{X}}$ is a minimal rational surface.    
\end{Prop}
\begin{pf}Rationality first.
Let $\tilde{X}$ be the compactification of ${X}$ after 
blowing down all the curves in
$\bigcup_{i\neq 0}D_{i}$, and 
still denote by  $D$ the image of $D_0$ after the blowing down process. 
Let $K$ denote the canonical divisor class of $\tilde{{X}}$. 
Since $D$ is a smooth rational curve we have,
from the genus formula (cf. \cite{GH})
$$
D\cdot(D+K)=2g(D)-2=-2
$$
and, since $D\cdot D=2$, it follows that
\begin{equation}\label{HK}
D\cdot K =-4.
\end{equation}
This implies that for all $m>1$ the line bundle
$\otimes^m K$ have no holomorphic sections, because any such 
would define an effective divisor which would have a 
non-negative intersection with $D$, contradicting  (\ref{HK}).
In particular,  
\begin{equation}\label{gg}
P_2:=\textit{ dim } H^0(\tilde{X}, \otimes^2 K)=0.\end{equation}
Since the first Betti number of $\tilde{X}$ is zero as well,
the rationality of $\tilde{{X}}$ follows from the rationality
criterion of Castelnuovo.

We now show that $\tilde{{X}}$
 satisfies
the ``minimality condition''.
This can be 
determined from the  "even parity" of the intersection form.
Indeed, according to a  well-known theorem of H. Wu, 
since $\tilde{X}$ is simply connected, 
the stronger requirement  that the square of any element in 
$H^2(\tilde{{X}}, {\mathbf Z})$ be even, is equivalent to 
$\tilde{{X}}$ being spin, i.e., 
its second Stiefel-Whitney class 
\begin{equation}\label{c12}
\omega_2(\tilde{X})=c_1(\tilde{{X}}) \mod{ (2) }=0.
\end{equation}
To show that indeed (\ref{c12})
holds in our case, recall that $c_1(\tilde{{X}})=-c_1(K)$.
Moreover, $H^2(\tilde{{X}})$ 
is generated over ${\mathbf{Q}}$ by $D$ and the class 
$S$ corresponding to $\tau^{-1}(0)$
and we may write
$$K=aD+dS, \quad a,b \in {\mathbf Q}.$$
Note also that $S\cdot S=-2$ (recall, the Euler number of the unit 
sphere bundle of the tangent bundle of $\tau^{-1}(0)$).

On the other hand, from Proposition \ref{euler} and Noether's Theorem 
(cf. \cite{GH}) which reads
$$1=\chi(\tilde{{X}}, {\cal O}_{\tilde{{X}}})=
\frac{1}{12}(K\cdot K +\epsilon (\tilde{{X}})),$$
we get $K\cdot K =8$,
hence $2a^2-2b^2=8$. But, since $\tau^{-1}(0)$ and the 
divisor at $\infty$ are disjoint, 
and hence their classes orthogonal, using (\ref{HK}), 
$$-4= D\cdot K = D\cdot (a D +b S)=2a.$$ 
Thus, we have $b=0$ and 
\begin{equation}\label{k-2h}
K=-2D,
\end{equation}
from which (\ref{c12}) follows, thus the  minimality
of $\tilde{{X}}$.
\end{pf}
\begin{Prop}
$\tilde{{X}}\simeq {\mathbf P}^1\times{\mathbf P}^1$.
\end{Prop}    
\begin{pf}    
It is a classical result in the theory of algebraic surfaces that 
a minimal rational surface, such as $\tilde{X}$,  
is 
either 
${\mathbf P}^2$ or belongs to a certain family of 
${\mathbf P}^1$ bundles ${\mathbf F}_k$ over ${\mathbf P}^1$, 
sometimes referred to 
as Hirzebruch surfaces (cf. \cite{GH}).
Clearly,
the first possibility is ruled out
because $\tilde{{X}}\setminus D$ is homotopically equivalent to 
$S^2$, in particular, non-contractible,
and thus we are left with the second possibility. 

Each ${\mathbf F}_k$, $0 \leq k$, 
is the projectivization of a rank 2 bundle, namely, 
$${\mathbf F}_k={\mathbf P}_{{\mathbf P}^1}({\cal O}_{{\mathbf P}^1}\oplus
(\otimes ^k{\cal O}_{{\mathbf P}^1}))
$$
where ${\cal O}_{{\mathbf P}^1}$ is the degree 1 bundle over 
${\mathbf P}^1$.
Given such a surface, the value of $k$ is determined by the  
self-intersection of the curve $C$ which is the image of the 
section defined by the subbundle 
$${\cal O}_{{\mathbf P}^1}\oplus\{0\} \subset 
{\cal O}_{{\mathbf P}^1}\oplus
(\otimes ^k{\cal O}_{{\mathbf P}^1})).$$
All we need is the fact that for any such surface, 
\begin{equation}\label{cc}
    C\cdot C\leq 0,
\end{equation}
    with equality 
if and only if the surface is ${\mathbf F}_0\simeq {\mathbf 
P}^1\times
{\mathbf P}^1$.

But, in our case we must have 
\begin{equation}\label{ch}
C\cdot D\geq 1
\end{equation}
since otherwise $C\cdot D =0$ and 
$C$ would be contained in the affine 
part ${X}= \tilde{{X}}\setminus D$, which is impossible.
It follows using  (\ref{cc}), (\ref{ch}) and 
(\ref{k-2h}) in the genus formula for $C$ that
$$0\leq C\cdot D-1=
\frac{1}{2} C\cdot C\leq 0,$$
and thus $C\cdot C=0$, which implies 
$$\tilde{X}\simeq F_{0}\simeq {\mathbf P}^{1}\times {\mathbf P}^{1}.$$
It follows that the $\sigma$ is given by $\sigma(z, w) = (\bar{w},
\bar{z})$, in suitable coordinates, since otherwise the fixed points
of $\sigma$ would be homeomorphic to the torus of dimension two. The
curve $D$ at infinity must be homologous to $C_1 + C_2$, where
the $C_i$ are the two factors in $\tilde{X}$. It follows that
$\tilde{X}$ is embedded as the usual quadric in ${\Bbb P}^3$, and that
$\sigma$ is conjugate linear in the homogeneous coordinates of ${\Bbb
P}^3$, with fixed point set homeomorphic to $S^2$. Without loss of
generality, we can conjugate $\tilde{X}$ by an element $g \in
{\mbox{Aut}}_{\Bbb C}(\tilde{X})$ so that $D$ is the standard diagonal
in $\tilde{X}$, the curve at infinity for the Grauert tube for the
round metric on $S^2$. The subgroup $G_D$ of  ${\mbox{Aut}}_{\Bbb
C}(\tilde{X})$ sending this $D$ into itself is isomorphic to the
adjoint group of $SO(3, {\Bbb C})$, and $G_D$ has two real forms
corresponding to $\sigma$ and the standard conjugation $\sigma_0$ of
$\tilde{X}$ coming from the round metric. The topology of the two sets
of real points in $\tilde{X}$ shows that both of the real forms of
$G_D$ are isomorphic to the real adjoint group of $SO(3, 1)$. Hence,
there is an element $g' \in G_D$ conjugating $\sigma$ to $\sigma_0$.
Hence, $X$ is biholomorphic to the standard affine quadric, and in
such a way that the complex conjugation $\sigma$ is equal to the
standard conjugation $\sigma_0$.

\end{pf}

%%%%%%%%%%%%%%%%%%%%%%%%%%%%%%%%%%%%%%%%%%%%%%%%%%%%%%%%%%%%%%%%%%%%%%%%%%%

\section{\bf{Open Questions}}

%%%%%%%%%%%%%%%%%%%%%%%%%%%%%%%%%%%%%%%%%%%%%%%%%%%%%%%%%%%%%%%%%%%%%%%%%%%

We have shown that Riemannian metrics on compact manifolds $M$ which give
rise to entire Grauert tubes $X$ must necessarily be algebraic, and the
Riemannian metric must extend to a rational metric on $\tilde{X}$,
regular and non-degenerate on $X$. In the previous section we showed
that this suffices in the lowest dimensional case to rigidify the
complex structure on $X$ for $M = S^2$. It remains to see why examples
of such metrics are so rare. 
It is a theorem of Sz\H{o}ke's
\cite{RSZOKE91} that among the surfaces of revolution on $S^2$ with
respect to a given rotation, the only ones corresponding to 
unbounded tubes form  a two-parameter family, where one of the parameters 
is simply rescaling the
metric. (See \cite{RMA4} for all known higher dimensional examples.) 
There are still quite a few rational metrics of the sort given
by theorem 1 and its refinements, in part because the tangent bundle
of our affine variety $X$ has a very large algebraic gauge group. It
is very difficult to encode the real information of $\tau$ and the
associated K\"ahler metric on $X$ in the metric $\tilde{g}$ on $X$.
It would be interesting even to see what the possibilities are on the
two sphere $S^2$.

It is not clear that there is enough known about minimal models in
three dimensions to be able to prove the analogue of theorem 2 for
metrics on the sphere $S^3$. 

Finally, it is an interesting question whether Demailly's method
\cite{DEMAILLY} can be used to resolve the algebraicization
conjecture of \cite{DB82}. The current paper settles this when $X$ is
a Grauert tube.

\end{document}